\documentclass[12pt]{amsart}
\usepackage{amssymb}
\usepackage{amscd}

\usepackage{amsmath}

\usepackage{amsfonts}

\newtheorem{lemma}{Lemma}[section]
\newtheorem{prop}[lemma]{Proposition}
\newtheorem{thm}[lemma]{Theorem}
\newtheorem{cor}[lemma]{Corollary}
\newtheorem{aplemma}{Lemma~A.\hspace{-1.5mm}}
\newtheorem{approp}{Proposition~A.\hspace{-1.5mm}}
\newtheorem{apthm}{Theorem~A.\hspace{-1.5mm}}
\newtheorem{apcor}{Corollary~A.\hspace{-1.5mm}}
\newtheorem{intthm}{Theorem}

\newtheorem{conj}[lemma]{Conjecture}

\newtheorem{rema}{Remark}

\newtheorem{defi}[lemma]{Definition}
\newtheorem{exa}[lemma]{Example}
\newtheorem{aprem}{Remark~A.\hspace{-1.5mm}}
\newtheorem{apdefi}{Definition~A.\hspace{-1.5mm}}
\newcommand{\bde}{\begin{defi}}
\newcommand{\ede}{\end{defi}\vspace{1mm}}
\newcommand{\ble}{\begin{lemma}}
\newcommand{\ele}{\end{lemma}}
\newcommand{\bpr}{\begin{prop}}
\newcommand{\epr}{\end{prop}}
\newcommand{\bt}{\begin{thm}}
\newcommand{\et}{\end{thm}}
\newcommand{\bco}{\begin{cor}}
\newcommand{\eco}{\end{cor}}
\newcommand{\bre}{\begin{rema}}
\newcommand{\ere}{\end{rema}}
\newcommand{\brea}{\begin{rema}}
\newcommand{\erea}{\end{rema}\vspace{1mm}}
\newcommand{\breb}{\begin{remb}}
\newcommand{\ereb}{\end{remb}\vspace{1mm}}
\newcommand{\bex}{\begin{exa}}
\newcommand{\eex}{\end{exa}}
\newcommand{\bpf}{\begin{proof}}
\newcommand{\epf}{\end{proof}\vspace{1mm}}

\newcommand{\bade}{\begin{apdefi}}
\newcommand{\eade}{\end{apdefi}}
\newcommand{\bale}{\begin{aplemma}}
\newcommand{\eale}{\end{aplemma}}
\newcommand{\bapr}{\begin{approp}}
\newcommand{\eapr}{\end{approp}}
\newcommand{\bat}{\begin{apthm}}
\newcommand{\eat}{\end{apthm}}
\newcommand{\baco}{\begin{apcor}}
\newcommand{\eaco}{\end{apcor}}
\newcommand{\bare}{\begin{aprem}}
\newcommand{\eare}{\end{aprem}}

\newcommand{\be}{\begin{enumerate}}
\newcommand{\ee}{\end{enumerate}}
\newcommand{\bcd}{\[\begin{CD}}
\newcommand{\ecd}{\end{CD}\]}
\newcommand{\bit}{\begin{itemize}}
\newcommand{\eit}{\end{itemize}}
\newcommand{\bq}{\begin{quote}}
\newcommand{\eq}{\end{quote}}
\newcommand{\ba}{\begin{array}}
\newcommand{\ea}{\end{array}}
\newcommand{\mcA}{\mathcal{A}}

\newcommand{\mcE}{\mathcal{E}}
\newcommand{\mcF}{\mathcal{F}}
\newcommand{\mcG}{\mathcal{G}}
\newcommand{\mcH}{\mathcal{H}}

\newcommand{\mcL}{\mathcal{L}}
\newcommand{\mcM}{\mathcal{M}}
\newcommand{\mcN}{\mathcal{N}}
\newcommand{\mcO}{\mathcal{O}}
\newcommand{\mcP}{\mathcal{P}}
\newcommand{\mcQ}{\mathcal{Q}}

\newcommand{\mcS}{\mathcal{S}}
\newcommand{\mcT}{\mathcal{T}}
\newcommand{\mcU}{\mathcal{U}}
\newcommand{\mcV}{\mathcal{V}}

\newcommand{\mbC}{\mathbb{C}}

\newcommand{\mbF}{\mathbb{F}}

\newcommand{\mbL}{\mathbb{L}}

\newcommand{\mbP}{\mathbb{P}}
\newcommand{\mbQ}{\mathbb{Q}}
\newcommand{\mbR}{\mathbb{R}}

\newcommand{\mbZ}{\mathbb{Z}}

\newcommand{\mfg}{\mathfrak{g}}

\newcommand{\migi}{\rightarrow}
\newcommand{\longmigi}{\longrightarrow}

\newcommand{\isom}{\stackrel{\sim}{\migi}}

\newcommand{\migiincl}{\hookrightarrow}

\newcommand{\migisurj}{\twoheadrightarrow}

 \newcommand{\XX}{X}

\newcommand{\Mgp}{\mcM_{g, \mbF_p}}

\newcommand{\N}[2]{\mcN_{#1}^{#2}}

\newcommand{\XY}{X}

\newcommand{\pic}[1]{\mr{Pic}^0_{#1}}
\newcommand{\ver}[1]{\mr{Ver}_{#1}}

\newcommand{\Dp}{\mcM^{^\text{Zzz...}}}

\newcommand{\triv}{\mcO}
\newcommand{\ind}{\circledast}

\newcommand{\mr}{\mathrm}
\newcommand{\hidden}[1]{\,}

\begin{document}

\title[the generic number of dormant indigenous bundles]{An explicit formula \\ for the generic number \\ of dormant indigenous bundles}
\author{Yasuhiro Wakabayashi}
\date{}
\maketitle
\footnotetext{Y. Wakabayashi: Research Institute for Mathematical Sciences, Kyoto University, Kyoto 606-8502, Japan;}
\footnotetext{e-mail: wakabaya@kurims.kyoto-u.ac.jp;}

\begin{abstract}

A dormant indigenous bundle is 
an  integrable $\mbP^1$-bundle on a proper hyperbolic curve of positive characteristic satisfying certain conditions.
Dormant indigenous bundles were
introduced and studied 
in  the $p$-adic Teichm\"{u}ller theory developed by  S. Mochizuki.
Kirti Joshi proposed a conjecture concerning an explicit  formula for  the  degree  over the moduli stack of curves of the moduli stack classifying dormant indigenous bundles.
In this paper,
we  give a proof for this conjecture of Joshi.
\end{abstract}
\tableofcontents 

\section*{Introduction}

Let 
\[\Dp_{g, \mbF_p}\]
 be the moduli stack classifying proper smooth curves of genus $g>1$ over $\mbF_p :=\mbZ/p\mbZ$ together with a {\it dormant} indigenous bundle (cf. the notation ``Zzz...''!).
 It is known (cf. Theorem 3.3) that  $\Dp_{g, \mbF_p}$ is represented by a smooth, geometrically connected Deligne-Mumford stack over $\mbF_p$ of dimension $3g-3$.
 Moreover, if we denote by $\Mgp$
 the moduli stack classifying proper smooth curves of genus $g$ over $\mbF_p$,
 then the natural projection $\Dp_{g, \mbF_p} \migi \Mgp$
 is finite, faithfully flat, and generically \'{e}tale.
 The main theorem of the present paper, which was conjectured by Kirti Joshi, asserts that if $p > 2(g-1)$, then the degree
$\mr{deg}_{\Mgp}{(\Dp_{g, \mbF_p})}$  of $\Dp_{g, \mbF_p}$ over $\Mgp$ may be calculated as follows:
\begin{intthm} [= Corollary 5.4]
 \[ \mr{deg}_{\Mgp}(\Dp_{g, \mbF_p})
 = \frac{ p^{g-1}}{2^{2g-1}} \cdot  \sum_{\theta =1}^{p-1}\frac{1}{\mr{sin}^{2g-2}(\frac{\pi \cdot  \theta}{p})}.  \] 
\end{intthm}
\vspace{5mm}

Here,  recall 
 that an indigenous bundle on a  proper smooth curve $X$ is a  $\mbP^1$-bundle on $X$, together with a  connection, which satisfies certain properties (cf. Definition 2.1).
The notion of an indigenous bundle was originally introduced and studied by Gunning in the context of compact hyperbolic Riemann surfaces (cf. ~\cite{G2}, \S\,2, p.\,69).
 One may think of an indigenous bundle as an algebraic object encoding uniformization data for $X$.
 It may be interpreted as a projective structure, i.e., a maximal atlas covered by coordinate charts on $X$ such that the transition functions are expressed as M\"{o}bius transformations.
Also, various equivalent mathematical objects, including  certain kinds of differential operators (related to Schwarzian equations) of kernel functions, 
have been studied by many mathematicians.

In the present paper, we focus on indigenous bundles in {\it positive characteristic}.
Just as in the case of the theory over $\mbC$, one may define the notion of an indigenous bundle and the moduli space classifying indigenous bundles.
Various properties of such objects  were firstly discussed in the context of the $p$-adic Teichm\"{u}ller  theory developed  by S. Mochizuki (cf. ~\cite{Mzk1}, ~\cite{Mzk2}).
(In a different point of view, Y. Ihara developed,  in, e.g.,  ~\cite{Ihara1}, ~\cite{Ihara2},  a theory of Schwarzian  equations in arithmetic context.)
One of the key ingredients in the development of this theory is
the study of  the $p$-curvature of indigenous bundles in characteristic $p$.
Recall that the $p$-curvature of a connection may be thought of as the obstruction to the  compatibility of $p$-power structures that appear in certain associated spaces of infinitesimal (i.e., ``Lie'')  symmetries.
We say that an indigenous bundle is {\it dormant} (cf. Definition 3.1) if its $p$-curvature  vanishes identically.
This condition on an indigenous bundle implies, in particular, the existence of ``sufficiently many'' horizontal sections locally in the Zariski topology.
Moreover, a dormant indigenous bundle corresponds, in a certain sense, to a certain type of rank $2$ semistable bundle.  Such semistable bundles have been studied in a different context (cf. \S\,6.1).
This sort of phenomenon is peculiar to the theory of indigenous bundles in {\it positive characteristic}.

In this context, one natural question is the following:
\begin{quote}\textit{Can one calculate explicitly the number of dormant indigenous bundles on a general curve?} \end{quote}
Since (as discussed above) $\Dp_{g, \mbF_p }$ is finite, faithfully flat, and generically \'{e}tale over $\Mgp$, 
the task of resolving this question may be reduced to the explicit computation of  the degree $\mr{deg}_{\Mgp}(\Dp_{g, \mbF_p})$  of $\Dp_{g, \mbF_p}$ over $\Mgp$.

In the case of $g=2$,
S. Mochizuki (cf. ~\cite{Mzk2}, Chap.\,V, \S\,3.2, p.\,267,  Corollary 3.7), H. Lange-C. Pauly (cf. ~\cite{LP}, p.\,180, Theorem 2), and B. Osserman (cf. ~\cite{O2}, p.\,274, Theorem 1.2)
verified  (by applying different methods) the equality
\[   \mr{deg}_{\mcM_{2, \mbF_p}}(\Dp_{2, \mbF_p }) = \frac{1}{24}\cdot (p^3-p). \]
For arbitrary $g$, Kirti Joshi conjectured, with his amazing insight, an explicit description, as  asserted in Theorem A, of the value  $ \mr{deg}_{\Mgp}(\Dp_{g, \mbF_p})$. 
(In fact, Joshi has proposed,  in personal communication to the author, a somewhat more general conjecture.
In the present paper,  however, we shall restrict our attention to a certain special case of this more general conjecture.)
The goal of the present paper is to verify the case $r =2$ of this conjecture of Joshi.

Our discussion in the present paper  follows, to a substantial extent, the ideas discussed in ~\cite{JP}, as well as in personal communication to the author by Kirti Joshi.
Indeed, certain of the  results obtained in the present paper are mild generalizations of the results obtained in ~\cite{JP} concerning rank $2$ opers to the case of  {\it families of curves} over quite general base schemes.
(Such {\it relative} formulations are necessary in the theory of the present paper, in order to consider {\it deformations} of various types of data.)
For example, Lemma 4.1 in the present paper corresponds to ~\cite{JP}, p.\,10, Theorem 3.1.6 (or ~\cite{JRXY}, \S\,5.3, p.\,627; ~\cite{SUN}, \S\,2, p.\,430, Lemma 2.1);
Lemma 4.2 corresponds to ~\cite{JP}, p.\,20, Theorem 5.4.1;
and Proposition 4.3 corresponds to ~\cite{JP}, p.\,21, Proposition 5.4.2.
Moreover, the insight concerning the connection with the formula of Holla (cf. Theorem 5.1), which is a special case of the Vafa-Intriligator formula,  is due to Joshi.

On the other hand, the new ideas introduced in the present paper may be summarized as follows.
First, we verify the {\it vanishing of obstructions} to deformation to characteristic zero  of 
a certain Quot-scheme that is related to $\Dp_{g, \mbF_p}$ (cf. Proposition 4.3, Lemma 4.4, and the discussion in the proof of Theorem 5.2).
Then we   relate the value $ \mr{deg}_{\Mgp}(\Dp_{g, \mbF_p})$
to the degree of the result of base-changing this Quot-scheme to $\mbC$ by applying 
{the formula of Holla (cf. Theorem 5.1, the proof of Theorem 5.2) {\it directly}.

Finally,  F. Liu and B. Osserman  
have shown (cf. ~\cite{LO}, p.\,126, Theorem 2.1) that
the value $ \mr{deg}_{\Mgp}(\Dp_{g, \mbF_p})$
may expressed as  a polynomial with respect to  the characteristic of the base field.
This was done by applying
Ehrhart's theory concerning the cardinality of the set of lattice points inside a polytope.
In \S\,6, we shall discuss the relation between this result and the main theorem of the present paper. 

\vspace{5mm}
\hspace{-4mm}{\bf Acknowledgement} \leavevmode\\
 \ \ \ The author cannot  express enough his sincere and deep gratitude to Professors Shinichi Mochizuki and Kirti Joshi (and hyperbolic curves of positive  characteristic!) for their  helpful suggestions and heartfelt encouragements, as well as for formulating  Joshi's conjecture.
Without their philosophies and amazing insights, his study of mathematics would have remained ``{\it dormant}".

The author would also like to thank  Professors Yuichiro Hoshi,  Brian Osserman, and Go Yamashita  for their  helpful discussions and advices. 
The author was  supported by the Grant-in-Aid for Scientific Research (KAKENHI No. 24-5691) and the Grant-in-Aid for JSPS fellows.

Special thanks go to Mr. Katsurou Takahashi, the staff members at ``{\it CAFE PROVERBS [15:17]}"  in Kyoto, Japan, and the various individuals with whom the author became acquainted there.
The author deeply appreciates the relaxed and comfortable environment that they provided for writing the  present paper.

Finally, the author would like to thank the referee for  reading carefully his manuscript and giving him some comments and suggestions. 

\vspace{10mm}

\section{Preliminaries} \vspace{0mm}

\subsection{}

Throughout this paper, we fix  an {\it odd} prime number $p$.

\subsection{}

We shall denote by $(Set)$ the category of (small) sets.  If $S$ is a Deligne-Mumford stack, then we shall denote by $(Sch)_S$  the category of schemes over $S$.

\subsection{}

If $S$ is a scheme and $\mcF$ an $\mcO_S$-module, then we shall denote by $\mcF^\vee$ its dual sheaf, i.e., $\mcF^\vee := \mcH om_{\mcO_S} (\mcF, \mcO_S)$.  
If $f : T \migi S$ is a finite flat scheme over a connected scheme $S$, then we shall denote by $\mr{deg}_S(T)$ the degree of $T$ over $S$, i.e., the rank of locally free $\mcO_S$-module $f_*\mcO_T$.

\subsection{}

If $S$ is a scheme (or more generally, a Deligne-Mumford stack), then
we define 
a {\it curve} over $S$ to be a geometrically connected and flat (relative) scheme $f : X \migi S$ over $S$ of relative dimension $1$.
Denote by $\Omega_{X/S}$ the  sheaf  of 1-differentials of $X$ over $S$ and
 $\mcT_{X/S}$ the dual sheaf of $\Omega_{X/S}$ (i.e., the sheaf of derivations of $X$ over  $S$).
We shall say that a proper smooth curve $f :X \migi S$ over $S$ is  {\it of genus $g$}
if the direct image $f_*\Omega_{X/S}$ is locally free of constant rank $g$.

\subsection{}
Let $S$ be a scheme over a field $k$, $X$  a smooth scheme over $S$,
$G$ an algebraic group over $k$, and $\mfg$  the Lie algebra of $G$.
Suppose that  $\pi : \mcE \migi X$ is a $G$-torsor over $X$.
Then we may associate to $\pi$  a short exact sequence
\[  0 \migi \mr{ad}(\mcE) \migi  \widetilde{\mcT}_{\mcE/S} \stackrel{\alpha_\mcE}{\migi} \mcT_{X/S} \migi 0,\]  
where $\mr{ad}(\mcE) := \mcE \times^{G} \mfg $ denotes the adjoint bundle associated to the $G$-torsor $\mcE$, and $\widetilde{\mcT}_{\mcE/S}$ denotes the subsheaf of $G$-invariant sections $(\pi_*\mcT_{\mcE/S})^G$ of $\pi_*\mcT_{\mcE/S}$. 
An {\it $S$-connection} on $\mcE$ is a split injection $\nabla : \mcT_{X/S} \migi  \widetilde{\mcT}_{\mcE/S}$ of the above short exact sequence (i.e., $\alpha_\mcE \circ \nabla = \mr{id}$).  If $X$ is of relative dimension $1$ over $S$, then any such $S$-connection is necessarily {\it integrable}, i.e.,
compatible with the Lie bracket structures on $\mcT_{X/S}$ and $\widetilde{\mcT}_{\mcE/S}= (\pi_*\mcT_{\mcE/S})^G$.

Assume that $G$ is a closed subgroup of $\mr{GL}_n$ for $n \geq 1$.
Then the notion of an $S$-connection defined here may be identified with
the usual definition of an $S$-connection on the associated vector bundle $\mcE \times^{G} (\mcO_X^{\oplus n})$ (cf. ~\cite{K}, p.\,10, Lemma 2.2.3; ~\cite{Kal}, p.\,178, (1.0)).
 In this situation, we shall not distinguish between these definitions of a connection.

If $\mcV$ is a vector bundle on $X$  equipped with an $S$-connection on $\mcV$, then we denote by $\mcV^\nabla$ the sheaf of horizontal sections in $\mcV$ (i.e., the kernel of the $S$-connection $\mcV \migi \Omega_{X/S}\otimes \mcV$).
\subsection{}

Let $S$ be a scheme of characteristic $p$ (cf. \S\,1.1) and $f:X \migi S$ a scheme over $S$.
The  {\it Frobenius twist of $X$ over $S$} is  the base-change $X^{(1)}$ of the $S$-scheme $X$
via  the absolute Frobenius morphism $F_S : S \migi S$ of $S$.
Denote by $f^{(1)} : X^{(1)} \migi S$ the structure morphism of the Frobenius twist of $X$ over $S$. 
The {\it relative Frobenius morphism of $X$ over $S$} is  the unique morphism $F_{X/S} : X \migi X^{(1)}$ over $S$ that fits into a commutative diagram of the form 
\[ \begin{CD}
X @> F_{X/S} >> X^{(1)} @>>> X
\\
@V f VV @V f^{(1)} VV @V f VV
\\
S @> \mr{id} >> S @>>> S,
\end{CD} \]
where the upper (respectively, the lower) composite is the absolute Frobenius morphism of $X$ (respectively, $S$).
If $f :X \migi S$ is smooth, geometrically connected and of relative dimension $n$,
then the relative Frobenius morphism $F_{X/S} : X \migi X^{(1)}$ is finite and faithfully flat of degree $p^n$.
In particular, the $\mcO_{X^{(1)}}$-module $F_{X/S *}\mcO_X$ is locally free of rank $p^n$.

\vspace{5mm}
\section{Indigenous bundles} \vspace{3mm}

In this section, we recall the notion of an indigenous bundle on a curve.
Much of the content of this section is 
implicit
in ~\cite{Mzk1}.


First, we discuss the definition of an indigenous bundle on a curve  (cf. ~\cite{Fr}, \S4, p.\,104; ~\cite{Mzk1}, Chap.\,I, \S\,2, p.\,1002, Definition 2.2).
Fix a scheme $S$ of characteristic $p$ (cf. \S 1.1)
and  a proper smooth curve $f : X \migi S$ of genus $g>1$ (cf. \S\,1.2).

\vspace{3mm}
\bde \leavevmode\\
\vspace{-5mm}
\begin{itemize}
\item[(i)]
Let $\mcP^\ind =(\mcP, \nabla)$ be a  pair consisting of a $\mr{PGL}_2$-torsor $\mcP$ over $X$ and  an (integrable) $S$-connection $\nabla$ on $\mcP$.
We shall say that $\mcP^\ind$ is an  {\it indigenous bundle} on $X/S$ if 
 there exists
 a globally defined section $\sigma$ of the associated $\mbP^1$-bundle
$\mbP^1_\mcP := \mcP \times^{\mr{PGL}_2} \mbP^1$
which has a nowhere vanishing derivative with respect to the connection $\nabla$.
We shall refer to the section $\sigma$ as the {\it Hodge section} of $\mcP^\ind$ (cf. Remark 2.1.1 (i)).
\item[(ii)]
Let $\mcP^\ind_1= (\mcP_1, \nabla_1)$, $\mcP^\ind_2=(\mcP_2, \nabla_2)$ be  indigenous bundles on $X/S$.
 An {\it isomorphism} from  $\mcP^\ind_1$ to $\mcP^\ind_2$ is an isomorphism $\mcP_1 \isom \mcP_2$ of $\mr{PGL}_2$-torsors  over $X$ that is compatible with the respective connections (cf. Remark 2.1.1 (iii)).
  \end{itemize}
  \ede

\begin{rema}  \leavevmode\\
Let $\mcP^\ind= (\mcP, \nabla)$ be an indigenous bundle on $X/S$. 
\begin{itemize}
\item[(i)]
The Hodge section $\sigma$
of $\mcP^\ind$
 is uniquely determined by the  condition that $\sigma$ have a nowhere vanishing derivative with respect to $\nabla$
 (cf. ~\cite{Mzk1}, Chap.\,I, \S\,2, p.\,1004, Proposition 2.4).
\item[(ii)]
The underlying $\mr{PGL}_2$-torsors of any two indigenous bundles on $X/S$ are isomorphic (cf. ~\cite{Mzk1}, Chap.\,I, \S\,2, p.\,1004, Proposition 2.5).
If there is a spin structure $\mbL = (\mcL, \eta_\mcL)$ on $X/S$ (cf. Definition 2.2),
then the $\mbP^1$-bundle $\mbP^1_\mcP$ is isomorphic to the projectivization  of an $\mbL$-bundle $\mcF$ as in Definition 2.3 (i), and the subbundle $\mcL \subseteq \mcF$ (cf. Definition 2.3 (i))
induces the Hodge section $\sigma$ (cf. Proposition 2.4).
\item[(iii)]
If  two indigenous bundles on $X/S$ are isomorphic, then any isomorphism between them is unique.
In particular, an indigenous bundle has no nontrivial automorphisms (cf. \S\,1.1;~\cite{Mzk1}, Chap.\,I, \S\,2, p.\,1006, Theorem 2.8).
\end{itemize}
 \end{rema}

\vspace{5mm}
 Next,  we  consider  a certain class of rank $2$ vector bundles
  with an integrable connection  (cf. Definition 2.3 (ii))
  associated to a specific choice of a spin structure (cf. Definition 2.2).
In particular, we show (cf. Proposition 2.4) that  such objects correspond to  indigenous bundles bijectively.
 We recall from, e.g., ~\cite{J}, \S\,2.1, p.\,25 the following:

\vspace{3mm}
\bde  \leavevmode\\
  A  {\it spin structure} on  $X/S$ is
 a pair 
 \[ \mbL := (\mcL, \eta_\mcL)\]
  consisting of  
 an invertible sheaf $\mcL$ on $X$ 
and  an isomorphism $\eta_\mcL : \Omega_{X/S} \isom \mcL^{\otimes 2}$.
A {\it spin curve} is a pair 
\[ (Y/S, \mbL)\]
consisting of a proper  smooth curve $Y/S$ of genus $g>1$ and a spin structure $\mbL$ on $Y/S$.
   \ede

\begin{rema}  \leavevmode\\
 \vspace{-5mm}
\begin{itemize}
\item[(i)]
$X/S$
 necessarily admits, at least \'{e}tale locally on $S$, a spin structure. 
Indeed, let us denote by $Pic_{\XX/S}^d$ the relative Picard scheme of $\XX/S$ classifying the set of (equivalence classes, relative to the equivalence relation determined by tensoring with a line bundle pulled back from the base $S$, of) degree $d$ invertible sheaves on $\XX$.
Then the  morphism 
\[  \mr{Pic}_{\XX/S}^{g-1} \migi \mr{Pic}_{\XX/S}^{2g-2} : [\mcL] \mapsto [\mcL^{\otimes 2}] \]
given by multiplication by $2$ is finite and \'{e}tale (cf. \S\,1.1).
Thus, the $S$-rational point of $\mr{Pic}_{\XX/S}^{2g-2}$ classifying the equivalence class $[\Omega^{}_{\XX/S}]$ determined by $\Omega^{}_{\XX/S}$  lifts,  \'{e}tale locally, to a point of $\mr{Pic}_{\XX/S}^{g-1}$.
\item[(ii)]
Let $\mbL = (\mcL, \eta_\mcL)$ be a spin structure on $\XX/S$ and
$T$ an $S$-scheme.
Then by  pulling back the structures $\mcL$, $\eta_\mcL$ via the natural projection $\XX \times_S T \migi \XX$, we obtain a spin structure on the curve $\XX \times_S T$ over $T$, which, by abuse of notation, we shall also denote by $\mbL$.   
\end{itemize}
\end{rema}

  In the following, let us fix a spin structure  $\mbL =(\mcL, \eta_\mcL)$ on $X/S$. 
  
\vspace{3mm}
\bde \leavevmode\\
\vspace{-5mm}
\begin{itemize}
\item[(i)]
An {\it $\mbL$-bundle} on $X/S$  is an extension, in the category of $\mcO_X$-modules,
\[ 0 \longmigi \mcL \longmigi \mcF \longmigi \mcL^\vee \longmigi 0 \]
of $\mcL^\vee$ by $\mcL$ 
whose restriction to each fiber over $S$  is nontrivial (cf. Remark 2.3.1 (i)).
We shall regard the underlying rank $2$ vector bundle associated to an $\mbL$-bundle as 
being equipped with a 2-step decreasing filtration $\{ \mcF^i \}_{i=0}^2$, namely, the filtration defined as follows:
\[   \mcF^2 :=0  \ \ \ \subseteq \ \ \ \mcF^1 := \mr{Im}(\mcL) \ \ \ \subseteq \ \ \  \mcF^0 := \mcF. \]
 \item[(ii)]
 An {\it $\mbL$-indigenous vector bundle} on $X/S$ is a triple
 \[  \mcF^\ind := (\mcF,  \nabla, \{ \mcF^1 \}_{i=0}^2) \]
 consisting of  an  $\mbL$-bundle $(\mcF, \{ \mcF^i \}_{i=0}^2)$ on $X/S$
 and an  $S$-connection $\nabla : \mcF \migi \mcF \otimes \Omega_{X/S} $ on $\mcF$ (cf. \S\,1.5)
satisfying the following two conditions.
\begin{itemize}
\item[(1)]
 If we equip $\mcO_X$ with the trivial connection and the determinant bundle $ \mr{det}(\mcF)$ with the natural connection induced by $\nabla$, then
 the natural composite isomorphism 
 \[ \mr{det}(\mcF) \isom \mcL \otimes \mcL^\vee \isom \mcO_X\]
  is horizontal.
   \item[(2)]
 The composite 
 \[ \mcL \stackrel{\nabla |_{\mcL}}{\migi} \mcF \otimes \Omega_{X/S} \migisurj  \mcL^\vee \otimes \Omega_{X/S}\]
of the restriction $\nabla |_{\mcL}$ of $\nabla$ to $\mcL  \ (\subseteq \mcF)$ and the morphism $\mcF \otimes \Omega_{X/S} \migisurj  \mcL^\vee \otimes \Omega_{X/S}$ induced by the quotient $\mcF \migisurj \mcL^\vee$
is an isomorphism.
This composite is often referred to as the {\it Kodaira-Spencer map}.
\end{itemize}
\item[(iii)]
Let $\mcF_1^\ind =(\mcF_1,  \nabla_1, \{ \mcF_1^1 \}_{i=0}^2)$, $\mcF_2^\ind =(\mcF_2,  \nabla_2, \{ \mcF_2^1 \}_{i=0}^2)$ on $X/S$ be $\mbL$-indigenous bundles on $X/S$.  Then an {\it isomorphism} from $\mcF_1^\ind$ to $\mcF_2^\ind$
is an isomorphism $\mcF_1 \isom \mcF_2$  of $\mcO_X$-modules that is compatible with the respective connections and filtrations and induces the identity morphism of $\mcO_X$ (relative to the respective natural composite isomorphisms discussed in (i)) upon taking determinants.
 \end{itemize}
  \ede

\begin{rema}  \leavevmode\\
 \vspace{-5mm}
 \begin{itemize}
 \item[(i)]
$\XY/S$ always admits an $\mbL$-bundle.  Moreover,
any two $\mbL$-bundles on $\XY/S$ are isomorphic Zariski locally on $S$.
Indeed, since $f:X\migi S$ is of relative dimension $1$,  the Leray-Serre spectral sequence 
$H^p(S, \mbR^qf_*\Omega_{X/S}) \Rightarrow H^{p+q}(X, f_*\Omega_{X/S})$
associated to the morphism $f:X\migi S$ yields an exact sequence
\[  0 \migi H^1(S, f_*\Omega_{X/S}) \migi  \mr{Ext}^1(\mcL^\vee, \mcL)  \migi H^0(S, \mbR^1f_*\Omega_{X/S})  \migi H^2(S, f_*\Omega_{X/S}) ,\] 
where the set  $\mr{Ext}^1(\mcL^\vee, \mcL)$ ($\cong H^1(X, \Omega_{X/S})$) corresponds to the set of  extension classes of $\mcL^\vee$ by $\mcL$.
In particular, 
if $S$ is an affine scheme, then the set of nontrivial extension classes corresponds bijectively to the set $H^0(S, \mcO_S)\setminus \{ 0\} \subseteq H^0(S, \mcO_S) \cong H^0(S, \mbR^1f_*\Omega_{X/S})$.

Also, we note that it follows immediately from the fact that the degree of the line bundle $\mcL$ on each fiber over $S$ is {\it positive} that the structure of $\mbL$-bundle on the underlying rank $2$ vector bundle of an $\mbL$-bundle is {\it unique}.

\item[(ii)]
If  two $\mbL$-indigenous vector bundles on $\XY/S$ are isomorphic, then any isomorphism between them is unique up to multiplication by an element of $\Gamma(S, \mcO_S)$ whose square is  equal to $1$ (i.e, $\pm 1$ if $S$ is connected).
In particular, the group of automorphisms of
an $\mbL$-indigenous vector bundle may be identified with the group of elements of $\Gamma(S,\mcO_S)$ whose square is equal to $1$.
(Indeed, these facts follow from  an argument  similar to the argument given in the proof in ~\cite{Mzk1}, Chap.\,I, \S\,2, p.\,1006, Theorem 2.8.)
\item[(iii)]
One may define, in an evident fashion, the pull-back of an $\mbL$-indigenous vector bundles on $X/S$ with respect to a morphism of schemes $S'\migi S$; this notion of pull-back is compatible, in the evident sense, with composites $S''\migi S'\migi S$.
\end{itemize}
\end{rema}

\vspace{5mm}
Let $\mcF^\ind =(\mcF, \nabla , \{ \mcF^i \}_{i=0}^2)$ be an $\mbL$-indigenous vector bundle on $X/S$.
By executing a change of structure group via the natural map $\mr{SL}_2\migi \mr{PGL}_2$, 
one may construct, from the pair $(\mcF,\nabla)$,
a $\mr{PGL}_2$-torsor $\mcP_\mcF$ together with an $S$-connection $\nabla_{\mcP_\mcF}$ on $\mcP_\mcF$.
Moreover, the subbundle $\mcL \ ( \subseteq \mcF)$ determines a globally defined section $\sigma$ of the associated $\mbP^1$-bundle $\mbP_\mcF^1:= \mcP_\mcF \times^{\mr{PGL}_2} \mbP^1$ on $X$.
One may verify easily from the condition  given in Definition 2.3 (ii) (2) that the pair $\mcP^\ind :=(\mcP_\mcF, \nabla_{\mcP_\mcF})$ forms an indigenous bundle on $X/S$, whose Hodge section is given by $\sigma$ (cf. Definition 2.1 (i)).
Then, we have (cf. ~\cite{Mzk1}, Chap.\,I, \S\,2, p.\,1004, Proposition 2.6) the following:
\vspace{3mm}
\bpr  \leavevmode\\
If $(X/S, \mbL)$ is a spin curve, then
the assignment $\mcF^\ind \mapsto \mcP^\ind$ discussed above determines 
a functor from the groupoid of $\mbL$-indigenous vector bundles on $X/S$ to the groupoid of indigenous bundles on $X/S$.
Moreover, this functor induces
a bijective correspondence between  the set of isomorphism classes of $\mbL$-indigenous vector bundles on $X/S$ (cf. Remark 2.3.1 (ii)) and the set of 
isomorphism classes of   indigenous bundles on $X/S$ (cf. Remark 2.1.1 (iii)).
 Finally, this correspondence is functorial with respect to $S$ (cf. Remark 2.3.1 (iii)). 
  \epr
\begin{proof}
The construction of a functor as asserted in the statement of Proposition 2.4 is routine.
The asserted (bijective) correspondence follows from ~\cite{Mzk1}, Chap.\,I, \S\,2, p.\,1004, Proposition 2.6.
 (Here, we note that Proposition 2.6 in loc.\,cit.\,states only that  an indigenous bundle determines an {\it indigenous vector bundle} (cf. ~\cite{Mzk1}, Chap.\,I, \S\,2, p.\,1002, Definition 2.2)  up to tensor product with a  line bundle together with a connection whose square is trivial.
But one may eliminate such an indeterminacy by
the condition that the underlying vector bundle be an $\mbL$-bundle.)
The asserted functoriality with respect to $S$ follows immediately from the construction of the assignment $\mcF^\ind \mapsto \mcP^\ind$ (cf. Remark 2.3.1 (iii)).
\end{proof}

\vspace{5mm}
\section{Dormant indigenous bundles} \vspace{3mm}

In this section, we recall
the notion of  a  dormant indigenous bundle
 and discuss various moduli functors related to this notion.

Let $S$ be a scheme over a field $k$ of characteristic $p$ (cf. \S\,1.1) and
$f : X \migi S$ a proper  smooth curve of genus $g>1$.
Denote by $X^{(1)}$ the Frobenius twist of $X$ over $S$  and
$F_{X/S} : X \migi X^{(1)}$ the relative Frobenius morphism of $X$ over $S$ (cf. \S\,1.6). 

First, we recall the definition of the $p$-curvature map.
Let us fix an algebraic group $G$ over $k$ and denote by $\mfg$ the Lie algebra of $G$.
Let $(\pi : \mcE \migi X,\nabla : \mcT_{X/S} \migi \widetilde{\mcT}_{\mcE/S})$ be a pair consisting of a $G$-torsor $\mcE$ over $X$ and an  $S$-connection $\nabla$ on $\mcE$, i.e., a section of the natural quotient $\alpha_\mcE : (\pi_*\mcT_{\mcE/S})^G=: \widetilde{\mcT}_{\mcE/S} \migi \mcT_{X/S}$ (cf. \S\,1.5).
If $\partial$ is a derivation corresponding to a local section $\partial$ of $\mcT_{X/S}$ (respectively, $\widetilde{\mcT}_{\mcE}:= (\pi_*\mcT_{\mcE/S})^G)$,
then we shall denote by $\partial^{[p]}$ the $p$-th iterate of $\partial$, which is also a derivation corresponding to a local section of $\mcT_{X/S}$ (respectively, $\widetilde{\mcT}_{\mcE}$).
Since $\alpha_\mcE(\partial^{[p]})= (\alpha_\mcE(\partial))^{[p]}$ for any local section of $\mcT_{X/S}$,
the image of the $p$-linear map from $\mcT_{X/S}$ to $\widetilde{\mcT}_{\mcE/S}$ defined by assigning $\partial \mapsto \nabla (\partial^{[p]})-(\nabla(\partial))^{[p]}$ is contained  in $\mr{ad}(\mcE)$ ($=\mr{ker}(\alpha_\mcE))$.
Thus, we obtain an $\mcO_X$-linear morphism
\[  \psi_{(\mcE, \nabla)} : \mcT_{X/S}^{\otimes p} \migi \mr{ad}(\mcE)  \]
determined by assigning
\[  \partial^{\otimes p} \mapsto  \nabla (\partial^{[p]})-(\nabla(\partial))^{[p]}. \]
We shall refer to the morphism $\psi_{(\mcE, \nabla)}$ as the {\it $p$-curvature map} of $(\mcE,\nabla)$.

If $\mcU$ is  a vector bundle on $X^{(1)}$,
then we may define an  $S$-connection (cf. \S\,1.5; ~\cite{Kal}, p.\,178, (1.0))
\[  \nabla^{\mr{can}}_\mcU :  F_{X/S}^*\mcU \migi F_{X/S}^*\mcU \otimes \Omega_{X/S}\] 
on the pull-back $F_{X/S}^*\mcU$ of $\mcU$,
which is uniquely determined by the condition that 
the sections of the subsheaf $F^{-1}_{X/S}(\mcU)$ be horizontal.
It is easily verified that the $p$-curvature map of  $(F_{X/S}^*\mcU, \nabla_\mcU^{\mr{can}})$ vanishes identically on $X$ (cf. Remark 3.0.1 (i)).

\begin{rema}  \leavevmode\\
Assume that $G$ is a closed subgroup of  $\mr{GL}_n$ for $n \geq1$ (cf. \S\,1.5).
Let $(\mcE, \nabla)$ be a pair consisting of a $G$-torsor $\mcE$ over $X$ and an $S$-connection $\nabla$ on $\mcE$.
Write $\mcV$ for the vector bundle on $X$ associated to $\mcE$ and 
 $\nabla_\mcV$ for the $S$-connection on $\mcV$ induced by $\nabla$.
\begin{itemize}
\item[(i)]
The $p$-curvature map $\psi_{(\mcE, \nabla)}$ of $(\mcE, \nabla)$
is compatible, in the evident sense, with the classical $p$-curvature map (cf., e.g., ~\cite{Kal}, \S\,5, p.\,190) 
of $(\mcV, \nabla_\mcV)$.
In this situation, we shall not distinguish between these definitions of the $p$-curvature map. 
\item[(ii)]
The sheaf $\mcV^\nabla$ of horizontal sections 
 in $\mcV$  may be considered as an $\mcO_{X^{(1)}}$-module via the underlying homeomorphism  of the relative Frobenius morphism $F_{X/S} : X \migi X^{(1)}$.
Thus, we have 
a natural horizontal morphism 
\[ \nu_{(\mcV,\nabla_\mcV)} : (F_{X/S}^*\mcV^\nabla, \nabla_{\mcV^\nabla}^{\mr{can}}) \migi (\mcV, \nabla_\mcV)\]
of $\mcO_{X}$-modules.
It is known (cf. ~\cite{Kal}, \S\,5, p.\,190, Theorem 5.1) that the $p$-curvature map of  $(\mcV, \nabla_\mcV)$ vanishes identically on $X$ if and only if $\nu_{(\mcV,\nabla_{\mcV})}$ is an isomorphism.
In particular,  the assignment $\mcV \mapsto (F^*_{X/S}\mcV, \nabla^{\mr{can}}_{\mcV^\nabla})$ determines an equivalence, which is compatible with the formation of tensor products (hence also symmetric and exterior products), between the category of vector bundles on $X^{(1)}$ and the category of vector bundles on $X$ equipped with an $S$-connection whose $p$-curvature  vanishes identically.
\end{itemize}
\end{rema}

\vspace{3mm}
\bde \leavevmode\\
 We shall say that an indigenous bundle $\mcP^\ind = (\mcP,\nabla)$ (respectively, an $\mbL$-indigenous vector bundle $\mcF^\ind =(\mcF, \nabla, \{ \mcF^i \}_{i=0}^2)$) on $X/S$
  is {\it dormant} if the $p$-curvature map of $(\mcP,\nabla)$ (respectively, $(\mcF, \nabla)$) vanishes identically on $X$.
 \ede

\vspace{3mm}
Next, we shall define a certain class of dormant indigenous bundles, which we shall refer to as {\it dormant ordinary}.
Let $\mcP^\ind = (\mcP, \nabla)$ be a dormant  indigenous bundle on $X/S$.
Denote by
\[  \mr{ad}(\mcP^\ind) := (\mr{ad}(\mcP), \nabla_{\mr{ad}}) \]
 the pair consisting of the adjoint bundle $\mr{ad}(\mcP)$
  associated to $\mcP$
  and the $S$-connection $\nabla_{\mr{ad}}$ on $\mr{ad}(\mcP)$
   naturally induced by
$\nabla$. 
Let us consider 
 the $1$-st relative de Rham cohomology sheaf 
$\mcH_{\mr{dR}}^1(\mr{ad}(\mcP^\ind))$
 of $\mr{ad}(\mcP^\ind)$,
 that is,
\[ \mcH_{\mr{dR}}^1(\mr{ad}(\mcP^\ind)) := \mbR^1f_*(\mr{ad}(\mcP)\otimes \Omega^\bullet_{X/S}), \]
where  $\mr{ad}(\mcP) \otimes \Omega^\bullet_{X/S}$
denotes the complex
 \[ \cdots \longmigi 0 \longmigi \mr{ad}(\mcP) \stackrel{\nabla_{\mr{ad}}}{\longmigi} \mr{ad}(\mcP) \otimes \Omega_{X/S} \longmigi 0 \longmigi \cdots \]
concentrated in degrees $0$ and $1$.
Recall (cf. ~\cite{Mzk1}, Chap.\,I, \S\,2, p.\,1006, Theorem 2.8) that there is a natural exact sequence
\[  0 \migi  f_*(\Omega_{X/S}^{\otimes 2}) \stackrel{}{\migi} \mcH_{\mr{dR}}^1(\mr{ad}(\mcP^\ind)) \stackrel{}{\migi} \mbR^1f_*(\mcT_{X/S}) \migi 0. \]
On the other hand, the natural inclusion
$\mr{ad}(\mcP)^\nabla \migiincl \mr{ad}(\mcP)$ of the subsheaf of horizontal sections induces a morphism of $\mcO_S$-modules
\[  \mbR ^1f_{*}(\mr{ad}(\mcP)^\nabla) \migi \mcH_{\mr{dR}}^1(\mr{ad}(\mcP^\ind)).  \]
Thus, by composing this morphism with the right-hand surjection in the above short exact sequence, we obtain a morphism
 \[ \gamma_{\mcP^\ind} : \mbR^1f_{*}(\mr{ad}(\mcP)^\nabla) \migi \mbR^1f_*(\mcT_{X/S})\] 
of $\mcO_S$-modules.
\vspace{3mm}
\bde \leavevmode\\
 \ \ \ We shall say that an indigenous bundle
$\mcP^\ind$ is 
  {\it dormant ordinary} if $\mcP^\ind$ is dormant and $\gamma_{\mcP^\ind}$ is an isomorphism.
   \ede

Next, let us introduce notations for various moduli functors classifying the objects discussed above.
Let $\Mgp$ be the moduli stack of proper smooth curves of genus $g>1$ {\it over $\mbF_p$}.
Denote by
\[   \mcS_{g, \mbF_p} : (Sch)_{\Mgp} \longmigi (Set)\]
(cf.  ~\cite{Mzk1}, Chap.\,I, \S\,3, p.\,1011,  the discussion preceding Lemma 3.2) the set-valued functor on $(Sch)_{\Mgp}$ (cf. \S\,1.2) which, to any  $\Mgp$-scheme $T$, classifying a curve $Y/T$, assigns the set of isomorphism classes  of indigenous bundles on $Y/T$.
Also, denote by
\[ \Dp_{g, \mbF_p}  \  \ \ (\text{resp.}, {^\circledcirc \Dp_{g, \mbF_p}}) \]
the subfunctor of $\mcS_{g, \mbF_p}$ classifying the set of isomorphism classes of dormant indigenous bundles (resp., dormant ordinary indigenous bundles).
By forgetting the datum of an indigenous bundle, we obtain natural transformations
\[   \mcS_{g, \mbF_p} \longmigi \Mgp, \ \ \ \ \  \Dp_{g, \mbF_p} \longmigi \Mgp. \]
 

Next, if $(X/S, \mbL)$ is a spin curve,
then we shall denote by
\[   \Dp_{X/S, \mbL} : (Sch)_S  \longmigi (Set) \] 
the set-valued functor on $(Sch)_S$ which, to any $S$-scheme $T$, assigns the set  of isomorphism classes of dormant $\mbL$-indigenous bundles on the curve $X \times_S T$ over $T$. 
It follows from Proposition 2.4 that there is a natural isomorphism of functors on $(Sch)_S$
\[  \Dp_{X/S, \mbL}  \isom \Dp_{g, \mbF_p} \times_{\Mgp} S, \]
where $\Dp_{g, \mbF_p} \times_{\Mgp} S$ denotes the fiber product of the natural projection  $\Dp_{g, \mbF_p} \migi \Mgp$ and the classifying morphism $S \migi \Mgp$ of $X/S$.

Next, we quote a result from $p$-adic Teichm\"uller theory due to S. Mochizuki concerning the moduli stacks (i.e., which are in fact {\it schemes}, relatively speaking, over $\Mgp$) that represent the functors discussed above.
Here, we wish to emphasize the importance of the {\it open density} of the dormant ordinary locus. 
As we shall see in Proposition 4.2 and its proof,  the properties stated in the following Theorem 3.3 enable us to relate a numerical calculation in {\it characteristic zero} to the degree of  certain moduli spaces of interest in {\it positive characteristic}.
\vspace{0mm}
\bt  \leavevmode\\
 \ \ \  
The functor $\mcS_{g, \mbF_p}$ is represented by  a relative affine space over $\Mgp$ of relative dimension $3g-3$.
 The functor $\Dp_{g, \mbF_p}$ is represented by a closed substack of $\mcS_{g, \mbF_p}$ which is finite and faithfully flat over $\Mgp$,
 and which is smooth and geometrically irreducible over $\mbF_p$.
 The functor ${^\circledcirc \Dp_{g, \mbF_p}}$ is an open dense substack of $\Dp_{g, \mbF_p}$ 
  and coincides with the \'{e}tale locus of
  $\Dp_{g,\mbF_p}$  over $\Mgp$.
  \et
\begin{proof}
The assertion follows from   ~\cite{Mzk1}, Chap.\,I, \S\,2, p.\,1007, Corollary 2.9;
 ~\cite{Mzk2}, Chap.\,II, \S\,2.3, p.\,152, Lemma 2.7;
 ~\cite{Mzk2}, Chap.\,II, \S\,2.3, p.\,153, Theorem 2.8 (and its proof). 
\end{proof}


In particular, 
it follows  that  it makes sense to speak of the {\it degree} 
\[ \mr{deg}_{\Mgp}(\Dp_{g, \mbF_p})\]
 of $\Dp_{g, \mbF_p}$
 over $\Mgp$.
The generic \'{e}taleness of $\Dp_{g, \mbF_p}$ over $\Mgp$ implies that
if $X$ is  a sufficiently generic proper smooth curve of genus $g$ over an algebraically closed field of characteristic $p$, then the number of dormant indigenous bundles on $X$ is exactly $\mr{deg}_{\Mgp}(\Dp_{g, \mbF_p})$.  
As we explained in the Introduction, our main interest in the present paper is  
the {\it explicit computation} of the value $\mr{deg}_{\Mgp}(\Dp_{g, \mbF_p})$.

\vspace{5mm}
\section{Quot-schemes} \vspace{3mm}

To calculate the value of $\mr{deg}_{\Mgp}(\Dp_{g, \mbF_p})$, it will be necessary to relate
 $\Dp_{g, \mbF_p}$ to certain 
Quot-schemes.
Here, to prepare for the discussion in \S 5 below, 
we introduce notions for Quot-schemes
in arbitrary characteristic.

Let $T$ be a noetherian scheme,
$Y$  a proper  smooth curve over $T$ of genus $g>1$ and $\mcE$ a vector bundle on $Y$. 
Denote by 
\[ \mcQ_{\mcE/Y/T}^{2, 0} : (Sch)_T \longmigi (Set)\]
the set-valued functor on $(Sch)_T$ which to any $f :T' \migi T$ associates the set of isomorphism classes of 
 injective morphisms of coherent $\mcO_{Y \times_T T'}$-modules
 \[  i :  \mcF \migi  \mcE_{T'},    \]
where $\mcE_{T'}$ denotes the pull-back of $\mcE$ via the projection $Y  \times_T T' \migi Y$, such that  the quotient $\mcE_{T'}/i (\mcF)$  is flat over $T'$ (which, since $Y/T$ is smooth of relative dimension $1$, implies that $\mcF$ is {\it locally free}), and $\mcF$ is of rank $2$ and degree $0$. 
It  is known (cf. ~\cite{FGA}, Chap.\,5, \S\,5.5, p.\,127, Theorem 5.14)  that $\mcQ_{\mcE/Y/T}^{2, 0}$ is represented by a proper  scheme over $T$.

Now let $(X/S, \mbL=(\mcL, \eta_\mcL))$ be a spin  curve of characteristic $p$
and denote, for simplicity, the relative Frobenius morphism $F_{X/S}: X \migi X^{(1)}$ by $F$.
Then in the following discussion, we consider the Quot-scheme discussed above
\[ \mcQ_{F_*(\mcL^\vee)/X^{(1)}/S}^{2, 0} \]
in the case where the data ``$(Y/T, \mcE)$" is taken to be $(X^{(1)}/S, F_*(\mcL^\vee))$.
If we denote by $\widetilde{i} : \widetilde{\mcF} \migi (F_*(\mcL^\vee))_{ \mcQ_{F_*(\mcL^\vee)/X^{(1)}/S}^{2, 0}}$ the tautological  injective morphism of sheaves on $X^{(1)} \times_S \mcQ_{F_*(\mcL^\vee)/X^{(1)}/S}^{2, 0}$, then the determinant bundle $\mr{det}(\widetilde{\mcF}) := \wedge^2(\widetilde{\mcF})$ determines a classifying morphism
\[  \mr{det} : \mcQ_{F_*(\mcL^\vee)/X^{(1)}/S}^{2, 0} \migi Pic_{X^{(1)}/S}^0 \]
to the relative Picard scheme $Pic_{X^{(1)}/S}^0$ (cf. Remark 2.2.1 (i)) classifying the set of equivalence classes of degree $0$ line bundles on  $X^{(1)}/S$.
We shall denote by
\[  \mcQ_{F_*(\mcL^\vee)/X^{(1)}/S}^{2, \triv} \]
 the scheme-theoretic inverse image, via $\mr{det}$, of the identity section of $Pic_{X^{(1)}/S}^0$.

Next, we discuss  a certain relationship between $\Dp_{X/S, \mbL}$ and $\mcQ_{F_*(\mcL^\vee)/X^{(1)}/S}^{2, \triv}$.
To this end, we introduce a certain filtered vector bundle with connection as follows.
Let us consider the rank $p$ vector bundle
\[ \mcA_\mcL := F^*F_*(\mcL^\vee)\]
on $X$ (cf. \S\,1.6), which has the canonical $S$-connection 
\[  \nabla^{\mr{can}}_{F_*(\mcL^\vee)}  \]
 (cf. the discussion preceding Remark 3.0.1).
By using this connection,  we may define a {\it $p$-step decreasing filtration} 
\[\{ \mcA_\mcL^i \}_{i=0}^p\]
 on
$\mcA_\mcL$
 as follows.
\begin{align*}
\mcA_\mcL^0  &:= \mcA_\mcL,\\
\mcA_\mcL^1 &:= \mr{ker}(\mcA_\mcL \stackrel{q}{\migisurj} \mcL^\vee), \\
\mcA_\mcL^j &:= \mr{ker}(\mcA_\mcL^{j-1} \stackrel{\nabla^{\mr{can}}_{F_*(\mcL^\vee)}|_{\mcA_\mcL^{j-1}}}{\longmigi} \mcA_\mcL \otimes \Omega_{X/S} \stackrel{}{\migisurj} \mcA_\mcL/\mcA_\mcL^{j-1} \otimes \Omega_{X/S})
\end{align*}
($j= 2, \cdots, p$),
where  $\mcA_\mcL (= F^*F_*(\mcL^\vee)) \stackrel{q}{\migisurj} \mcL^\vee$ denotes the  natural quotient determined by 
the adjunction relation ``$F^*(-) \dashv F_*(-)$"
(i.e., ``the functor $F^*(-)$ is left adjoint to the functor $F_*(-)$").

\vspace{3mm}
\ble 
 \leavevmode\\
 \vspace{-5mm}
 \begin{itemize}
 \item[(i)]
 For each  $j= 1, \cdots p-1$, the map
 \[\mcA_\mcL^{j-1}/ \mcA_\mcL^j \migi \mcA_\mcL^{j}/ \mcA_\mcL^{j+1} \otimes \Omega_{X/S}\]
  defined by assigning $\overline{a} \mapsto \overline{\nabla_{F_*(\mcL^\vee)}^{\mr{can}}(a)}$ ($a \in \mcA_\mcL^{j-1}$), where the 
``bars" denote the images in the respective quotients,
 is well-defined and determines an  isomorphism of $\mcO_X$-modules.
 \item[(ii)]
Let us  identify $\mcA_\mcL^1/\mcA_\mcL^2$ with $\mcL$ via the isomorphism  
 \[\mcA_\mcL^1/\mcA_\mcL^2 \isom \mcA_\mcL^0/\mcA_\mcL^1 \otimes \Omega_{X/S} \isom \mcL^\vee \otimes \Omega_{X/S} \isom\mcL,\]
obtained by composing the isomorphism of (i) (i.e., the first isomorphism of the display)
with the tautological isomorphism arising from the definition of $\mcA_\mcL^1$ (i.e., the second isomorphism of the display), followed by the isomorphism determined by the given spin structure (i.e., the third isomorphism of the display).
 Then the natural extension structure
 \[  0 \migi  \mcA_\mcL^1/\mcA_\mcL^2 \migi \mcA_\mcL/\mcA_\mcL^2 \migi \mcA_\mcL/\mcA_\mcL^1 \migi 0\]
determines a structure of $\mbL$-bundle on $\mcA_\mcL/\mcA_\mcL^2$.
 \end{itemize}
 \ele
\begin{proof}
The various assertions of Lemma 4.1 follow from an argument (in the case where $S$ is an arbitrary scheme) similar to
the argument (in the case where $S= \mr{Spec}(k)$ for an algebraically closed field $k$) given in the proofs of ~\cite{JRXY}, \S\,5.3, p.\,627 and ~\cite{SUN}, \S\,2, p.\,430, Lemma 2.1.
\end{proof}

\vspace{3mm}
\ble 
 \leavevmode\\
 \ \ \ Let $g : \mcV \migi F_*(\mcL^\vee)$ be an injective morphism classified by an $S$-rational point of $\mcQ_{F_*(\mcL^\vee)/X^{(1)}/S}^{2, 0}$  and denote by 
 $\{ (F^*\mcV)^i \}_{i=0}^p$ the  filtration on the pull-back $F^*\mcV$ defined by setting 
 \[ (F^*\mcV)^i := (F^*\mcV) \cap (F^*g)^{-1}(\mcA_\mcL^i), \]
where we denote by $F^*g$ the pull-back of $g$ via $F$.
\begin{itemize}
\item[(i)]
 The composite
 \[ F^*\mcV \migi \mcA_\mcL \migisurj \mcA_\mcL /\mcA_\mcL^2 \]
of $F^*g$ with the natural quotient $\mcA_\mcL \migisurj \mcA_\mcL /\mcA_\mcL^2 $  is an isomorphism of $\mcO_X$-modules.
\item[(ii)]
If, moreover, $g$ corresponds to an $S$-rational point of $\mcQ_{F_*(\mcL^\vee)/X^{(1)}/S}^{2, \triv}$, then the triple 
 \[(F^*\mcV, \nabla^{\mr{can}}_\mcV, \{ (F^*\mcV)^i\}_{i=0}^2),\]
where $\nabla^{\mr{can}}_\mcV$ denotes the canonical connection on $F^*\mcV$ (cf. the discussion preceding Remark 3.0.1),
  forms a dormant $\mbL$-indigenous bundle on $X/S$. 
\end{itemize}
 \ele
\begin{proof}
First, we consider  assertion (i).
Since $F^*\mcV$ and $\mcA_\mcL /\mcA_\mcL^2$ are flat over $S$, it suffices, by  
considering the various fibers
over $S$, to verify the case where $S = \mr{Spec}(k)$ for  a field $k$.  
If we write $\mr{gr}^i := (F^*\mcV)^i/(F^*\mcV)^{i+1}$ ($i= 0, \cdots, p-1$),
 then it follows immediately from the definitions that
 the coherent $\mcO_X$-module $\mr{gr}^i$ admits a natural
 embedding 
 \[  \mr{gr}^i \migiincl \mcA_\mcL^i/\mcA_\mcL^{i+1} \]
into the subquotient $\mcA_\mcL^i/\mcA_\mcL^{i+1}$.
 Since this subquotient is a line bundle (cf. Lemma 4.1 (i), (ii)), one verifies easily that
 $\mr{gr}^i$ is either trivial or a line bundle.
In particular, since $F^*\mcV$ is of rank $2$, the cardinality of the set $I := \big\{ i \big| \mr{gr}^i \neq 0\big\}$ is exactly $2$.
Next, let us observe that the pull-back $F^*g$ of $g$ via $F$ is compatible with the respective connections $\nabla^{\mr{can}}_\mcV$ (cf. the statement of assertion (ii)), $\nabla^{\mr{can}}_{F_*(\mcL^\vee)}$.  
Thus, it follows from Lemma 4.1 (i) that $\mr{gr}^{i+1} \neq 0$ implies $\mr{gr}^i \neq 0$.
But this implies that  $I= \{ 0,1\}$, and hence that
 the composite
\[  F^*\mcV \migi \mcA_\mcL \migisurj \mcA_\mcL /\mcA_\mcL^2 \]
is an isomorphism at the generic point of $X$.
On the other hand, observe that
\[ 
\mr{deg}(F^*\mcV) =  p \cdot \mr{deg}(\mcV) = p \cdot 0 = 0 
 \]
and 
\[ \mr{deg}(\mcA_\mcL/\mcA_\mcL^2)  = \mr{deg}(\mcA_\mcL/\mcA_\mcL^1) + \mr{deg}(\mcA_\mcL^1/\mcA_\mcL^2)  = \mr{deg}(\mcL^\vee) + \mr{deg}(\mcL) = 0\]
 (cf. Lemma 4.1 (i)).
Thus,  by comparing the  respective degrees of $F^*\mcV$ and $\mcA_\mcL/\mcA_\mcL^2$,
we conclude that the above composite is  an isomorphism of $\mcO_X$-modules.  This completes the proof of  assertion (i).
Assertion  (ii) follows immediately from the definition of an $\mbL$-indigenous bundle, assertion (i),  and Lemma 4.1 (i), (ii).
\end{proof}
\vspace{3mm}

By applying the above lemma, we may conclude that the moduli space $\Dp_{X/S, \mbL}$ is isomorphic to the Quot-scheme $\mcQ_{F_*(\mcL^\vee)/X^{(1)}/S}^{2, \triv}$ as follows.

\vspace{3mm}
\bpr 
\leavevmode\\
 \ \ \  
Let $(X/S, \mbL)$ be a spin curve.
Then
there is an isomorphism  of $S$-schemes
\[ \mcQ_{F_*(\mcL^\vee)/X^{(1)}/S}^{2, \triv}   \isom \Dp_{X/S, \mbL}.  \] 
 \epr
\begin{proof}
The assignment 
\[ [g : \mcV \migi F_*(\mcL^\vee)] \mapsto (F^*\mcV, \nabla^{\mr{can}}_{F^*\mcV}, \{ (F^*\mcV)^i\}_{i=0}^2),\]
discussed in Lemma 4.2,
determines  (by Lemma 4.2 (ii)) a map
\[ \alpha_S : \mcQ_{F_*(\mcL^\vee)/X^{(1)}/S}^{2, \triv}(S) \migi \Dp_{X/S, \mbL} (S)\]
between the respective sets of   $S$-rational points.
By the functoriality of the construction of $\alpha_S$ with respect to $S$,
it suffices to prove the bijectivity of   $\alpha_S$.

The {\it injectivity} of $\alpha_S$ follows from the observation that
any  element $[g : \mcV \migi F_*(\mcL^\vee)] \in \mcQ_{F_*(\mcL^\vee)/X^{(1)}/S}^{2, \triv}(S)$
is, by adjunction,  determined  by  the morphism $F^*\mcV \migi \mcL^\vee$, i.e.,
the natural surjection, as  in Definition 2.3 (i),  arising from the fact that  $F^*\mcV$ is an $\mbL$-bundle (cf. Lemma 4.2 (ii)).

Next, we consider the {\it surjectivity} of $\alpha_S$.
Let  $(\mcF, \nabla, \{ \mcF^i \}_i)$ be a dormant $\mbL$-indigenous bundle  on $X/S$.
Consider
the composite $F^*\mcF^\nabla \isom \mcF \migisurj \mcL^\vee$
of the natural horizontal isomorphism $F^*\mcF^\nabla \isom \mcF$ (cf. Remark 3.0.1 (ii)) with the natural surjection $\mcF \migisurj \mcF/\mcF^1 =\mcL^\vee$.  This composite determines
 a morphism 
 \[ g_\mcF : (\mcF \cong) F^*\mcF^\nabla \migi F^*F_*(\mcL^\vee) (=:\mcA_\mcL)\] 
via  the adjunction relation ``$F^*(-) \dashv F_*(-)$" (cf, the discussion preceding Lemma 4.1) and pull-back by $F$.

Next, we claim that $g_\mcF$ is {\it injective}.
Indeed, since $g_\mcF$ is (tautologically, by construction!) compatible with the respective surjections 
$\mcF \migisurj \mcL^\vee$, $\mcA_\mcL \migisurj \mcL^\vee$ to $\mcL^\vee$,
we conclude  that
$g_\mcF (\mcF^1) \subseteq \mcA_\mcL^1$, and $\mr{ker}(g_\mcF) \subseteq \mcF^1$.
Since $g_\mcF$ is {\it manifestly horizontal} (by construction), $\mr{ker}(g_\mcF)$ is stabilized by $\nabla$, hence contained in the kernel of the Kodaira-Spencer map $\mcF^1\migi \mcF/\mcF^1 \otimes \Omega_{X/S}$(cf. Definition 2.3 (ii) (2)), which is an isomorphism by the definition of an $\mbL$-indigenous bundle (cf. Definition 2.3 (ii)).
This implies that $g_\mcF$ is injective and completes the proof of the claim.
Moreover, by applying a similar  argument  to the pull-back of $g_\mcF$ via any base-change over $S$, one concludes that $g_\mcF$ is {\it universally injective} with respect to base-change over $S$.
This implies  that $\mcA_\mcL/g_\mcF(\mcF)$ is {\it flat} over $S$ (cf.~\cite{MAT}, p.\,17, Theorem 1).

Now denote by $g_\mcF^\nabla : \mcF^\nabla \migi F_*(\mcL^\vee)$
the morphism obtained by restricting $g_\mcF$ to the respective subsheaves of horizontal sections in $\mcF$, $\mcA_\mcL$.
Observe that the pull-back of $g_\mcF^\nabla$ via $F$ may be identified with $g_\mcF$, and that $F^*(F_*(\mcL^\vee) /g^\nabla_\mcF(\mcF^\nabla))$ is naturally isomorphic to $\mcA_\mcL/g_\mcF(\mcF)$.
Thus, it follows from the faithful flatness of $F$ that $g_\mcF^\nabla$ is injective, and $F_*(\mcL^\vee)/g^\nabla_\mcF(\mcF^\nabla)$ is flat over $S$.
On the other hand, since the determinant of $(\mcF, \nabla)$ is  trivial, $\mr{det}(\mcF^\nabla)$ is isomorphic to the trivial $\mcO_{X^{(1)}}$-module (cf. Remark 3.0.1 (ii)).
Thus,  $g_\mcF^\nabla$ determines an $S$-rational point of $\mcQ_{F_*(\mcL^\vee)/X^{(1)}/S}^{2, \triv}$ that  is mapped by $\alpha_S$ to the $S$-rational point of $\Dp_{X/S, \mbL}$ corresponding to $(\mcF, \nabla, \{ \mcF^i \}_i)$.
This implies that  $\alpha_S$ is surjective and  hence completes  the proof of Proposition 4.3. 
\end{proof}
\vspace{3mm}

Next, we relate $\mcQ_{F_*(\mcL^\vee)/X^{(1)}/S}^{2, \triv}$ to $\mcQ_{F_*(\mcL^\vee)/X^{(1)}/S}^{2, 0}$.
By pulling back line bundles on $X^{(1)}$
via the relative Frobenius $F_{} : X \migi X^{(1)}$, we obtain   a morphism 
\[ \pic{X^{(1)}/S} \migi \pic{X/S}  \]
\[   [\mcN] \mapsto [F_{}^*\mcN]. \]
We shall denote by
\[  \ver{X/S} \]
the scheme-theoretic inverse image, via this morphism,  of the identity section of $\pic{X/S}$.
It is well-known  (cf. ~\cite{DG}, EXPOSE VII, \S\,4.3, pp.\,440-443; ~\cite{MIL}, \S\,8, p.\,114, Proposition 8.1 and p.\,115, Theorem 8.2; ~\cite{MES}, APPENDIX, p.\,175, Lemma (1.0))
 that $\ver{X/S}$ is finite and faithfully flat over $S$ of degree $p^g$ and, moreover, \'{e}tale 
 over the points $s$ of $S$ such that the fiber of $X/S$ at $s$ is ordinary.
(Recall that the locus of $\Mgp$ classifying ordinary  curves is open and dense.)
Then we have the following

\vspace{3mm}
\ble \leavevmode\\
There is an isomorphism of $S$-schemes
\[   \mcQ_{F_*(\mcL^\vee)/X^{(1)}/S}^{2, \triv} \times_{S}   \ver{X/S} \isom \mcQ_{F_*(\mcL^\vee)/X^{(1)}/S}^{2, 0}. \]
 \ele
\begin{proof}
It  suffices to prove  that 
there is a bijection between 
the respective sets of $S$-rational points that is functorial with respect to $S$.

Let $(g :\mcV \migi F_*(\mcL^\vee), \mcN)$ be an element of $(\mcQ_{F_*(\mcL^\vee)/X^{(1)}/S}^{2, \triv} \times_{S}   \ver{X/S})(S)$.
It follows from the projection formula that the composite 
\[   g_\mcN : \mcV \otimes \mcN \migi F_*(\mcL^\vee) \otimes \mcN \migi F_*(\mcL^\vee \otimes F^*\mcN) \isom F_*(\mcL^\vee \otimes \mcO_X) = F_*(\mcL^\vee) \]
determines an element of $\mcQ_{F_*(\mcL^\vee)/X^{(1)}/S}^{2, 0}(S)$.
Thus, we obtain a functorial (with respect to $S$) map 
\[\gamma_S :(\mcQ_{F_*(\mcL^\vee)/X^{(1)}/S}^{2, \triv} \times_{S}   \ver{X/S})(S) \migi \mcQ_{F_*(\mcL^\vee)/X^{(1)}/S}^{2, 0}(S).\]
Conversely, let $g :\mcV \migi F_*(\mcL^\vee)$ be an injective morphism classified by an element of $\mcQ_{F_*(\mcL^\vee)/X^{(1)}/S}^{2, 0}(S)$.
Consider the injective morphism $g_{\mr{det}(\mcV)^{\otimes \frac{p-1}{2}}}$, i.e., the morphism $g_\mcN$ constructed above in the case where  ``$\mcN$'' is taken to be  $\mcN = \mr{det}(\mcV)^{\otimes \frac{p-1}{2}}$.
Here, we observe that
\[  \mr{det}(\mcV \otimes \mr{det}(\mcV)^{\otimes \frac{p-1}{2}}) \cong \mr{det}(\mcV) \otimes  \mr{det}(\mcV)^{\otimes 2 \cdot \frac{p-1}{2}} \cong \mr{det}(\mcV)^{\otimes p} \cong F_S^*(F^*(\mr{det}(\mcV))), \]
where $F_S^*(-)$ denotes the pull-back by the morphism $X^{(1)} \migi X$ obtained by base-change of $X/S$ via the absolute Frobenius morphism $F_S : S \migi S$ of $S$ (cf. \S\,1.6).
On the other hand, since $F^*(\mr{det}(\mcV)) \cong (\mcA_\mcL/\mcA_\mcL^1) \otimes (\mcA_\mcL^1/\mcA_\mcL^2) \cong \mcL^\vee \otimes \mcL \cong \mcO_X$ (cf. Lemmas 4.1 (ii), 4.2 (i)),  it follows that the determinant of 
$\mcV \otimes \mr{det}(\mcV)^{\otimes \frac{p-1}{2}}$ is trivial.
Thus the pair $(g_{\mr{det}(\mcV)^{\otimes \frac{p-1}{2}}}, \mr{det}(\mcV))$ determines an element of $(\mcQ_{F_*(\mcL^\vee)/X^{(1)}/S}^{2, \triv} \times_{S}   \ver{X/S})(S)$.
One verifies easily that
this assignment determines an inverse to $\gamma_S$.
This completes the proof of Lemma 4.4. 
\end{proof}
\vspace{3mm}

\vspace{6mm}
\section{Computation via the Vafa-Intriligator formula} \vspace{3mm}

By combining  Proposition 4.3, Lemma 4.4, and  the discussions preceding Theorem 3.3 and Lemma 4.4,
we obtain the following equalities:
\[  \mr{deg}_{\Mgp}(\Dp_{g, \mbF_p}) =  \mr{deg}_S(\Dp_{X/S, \mbL}) = \mr{deg}_S( \mcQ_{F_*(\mcL^\vee)/X^{(1)}/S}^{2, \triv})  = \frac{1}{p^g} \cdot \mr{deg}_S(\mcQ_{F_*(\mcL^\vee)/X^{(1)}/S}^{2, 0}). \]
Therefore,  to determine the value of  $ \mr{deg}_{\Mgp}(\Dp_{g, \mbF_p})$, it suffices to calculate the value $\mr{deg}_S(\mcQ_{F_*(\mcL^\vee)/X^{(1)}/S}^{2, 0}) $ (for an {\it arbitrary} spin  curve $(X/S, \mbL)$).

In this section, we review a numerical formula concerning  the degree of a certain Quot-scheme over the field of complex number $\mbC$ and relate it to the degree of the Quot-scheme in positive characteristic.

 Let $C$ be a smooth proper curve over $\mbC$ of genus $g>1$.
If $r$ is an integer, and $\mcE$ is
 a vector bundle  on $C$ of rank $n$ and degree $d$
with $1 \leq r \leq n$,
 then we define invariants
 \[ \begin{split} e_{\mr{max}}(\mcE,r) & := \mr{max} \big\{  \mr{deg}(\mcF) \in \mbZ\ \big|\ \text{$\mcF$ is a subbundle of $\mcE$ of rank $r$ }   \big\},   \\
  s_r(\mcE) & := d \cdot r - n \cdot e_{\mr{max}}(\mcE,r). \end{split} \]
(Here, we recall that one verifies immediately, for instance, by considering an embedding of $\mcE$ into a direct sum of $n$ line bundles, that $e_{\mr{max}}(\mcE,r) $ is well-defined.)

In the following, we review some facts concerning these invariants (cf. ~\cite{HIR}; ~\cite{LN}; \ ~\cite{H}).
Denote by ${^{\text{s}}\mcN}_C^{n,d}$ the moduli space of stable bundles on $C$ of rank $n$ and degree $d$ (cf. ~\cite{LN}, \S\,1, pp.\,310-311).
It is known that  ${^{\text{s}}\mcN}_C^{n,d}$ is irreducible  (cf.  the discussion at the beginning of ~\cite{LN}, \S\,2, p.\,311).
Thus, it makes sense to speak of a ``sufficiently general" stable bundle in ${^{\text{s}}\mcN}_C^{n,d}$, i.e., a stable bundle that corresponds to a point of the scheme ${^{\text{s}}\mcN}_C^{n,d}$ that lies outside some fixed closed subscheme.
If $\mcE$ is  a sufficiently general stable bundle in ${^{\text{s}}\mcN}_C^{n,d}$,
then it holds (cf. ~\cite{LN}, \S\,1, pp.\,310-311) that
$ s_r (\mcE) = r (n-r)(g-1) + \epsilon$, where
$\epsilon$ is the unique integer such that $ 0 \leq \epsilon < n$ and $s_r(\mcE) = r \cdot d$ mod $n$.
Also, 
the number $\epsilon$
coincides (cf. ~\cite{H}, \S\,1, pp.\,121-122) with the dimension of every irreducible component of the Quot-scheme 
$\mcQ^{r, e_\mr{max}(\mcE,r)}_{\mcE/C/\mbC}$ (cf. \S\,4).
If, moreover, the equality $s_r(\mcE) =r(n-r)(g-1)$ holds (i.e., $\mr{dim}(\mcQ^{r, e_\mr{max}(\mcE,r)}_{\mcE/C/\mbC})=0$), then $\mcQ^{r, e_\mr{max}(\mcE,r)}_{\mcE/C/\mbC}$ is \'{e}tale over $\mr{Spec}(\mbC)$ (cf. ~\cite{H}, \S\,1, pp.\,121-122).
Finally,  under this particular assumption, a formula for the degree of this Quot-scheme was given by Holla
as follows. 
\vspace{3mm}
\bt   \leavevmode\\
 \ \ \  Let $C$ be a  proper smooth curve over $\mbC$ of genus $g>1$, $\mcE$  a sufficiently general stable bundle in ${^{\text{s}}\mcN}_C^{n,d}$.
Write $(a,b)$ for the unique pair of integers such that $d =an-b$ with $0 \leq b < n$.
Also, we suppose that the equality $s_r(\mcE) = r(n-r)(g-1)$ (equivalently, $e_\mr{max}(\mcE,r) = (dr - r(n-r)(g-1))/n$) holds.
 Then we have 
 \[  \mr{deg}_{\mbC}(\mcQ^{r,e_\mr{max}(\mcE,r)}_{\mcE/C/\mbC}) = \frac{(-1)^{(r-1)(br-(g-1)r^2)/n}n^{r(g-1)}}{r!}\sum_{\rho_1 , \cdots , \rho_r}\frac{(\prod_{i=1}^r \rho_i)^{b-g+1}}{\prod_{i \neq j} (\rho_i - \rho_j)^{g-1}}, \]
 where $\rho_i^n =1$, for $1 \leq i \leq r$ and the sum is over tuples $(\rho_1 , \cdots , \rho_r)$ with $\rho_i \neq \rho_j$.
 \et
\begin{proof}
The assertion follows from  ~\cite{H}, \S\,4, p.\,132, Theorem 4.2, where ``$k$'' (respectively, ``$r$'') corresponds to our $r$ (respectively, $n$).
\end{proof}


By applying this formula,
we  conclude the same kind of formula 
for certain vector bundles 
in
positive characteristic,
as follows.
\vspace{3mm}
\bt \leavevmode\\
 \ \ \ 
 Let  $k$ an algebraically closed field of characteristic $p$ and
$(X/k, \mbL=(\mcL, \eta_\mcL))$ 
a  spin curve of genus $g>1$.
Suppose that $X/k$ is sufficiently general in $\Mgp$.
(Here, we recall that $\Mgp$ is irreducible (cf. ~\cite{DM}, \S\,5); thus, it makes sense to speak of a ``sufficiently general" $X/k$, i.e., an $X/k$ that determines a point of $\Mgp$ that lies outside some fixed closed substack.)
Then  $\mcQ_{F_*(\mcL^\vee)/X^{(1)}/k}^{2, 0}$ is finite and \'{e}tale  over $k$.
If, moreover,  we suppose that $p > 2(g-1)$, then
 the degree
 $\mr{deg}_k(\mcQ_{F_*(\mcL^\vee)/X^{(1)}/k}^{2, 0})$ of $\mcQ_{F_*(\mcL^\vee)/X^{(1)}/k}^{2, 0}$ over $\mr{Spec}(k)$ is given by the following formula:
\[ \begin{split} \mr{deg}_k(\mcQ_{F_*(\mcL^\vee)/X^{(1)}/k}^{2, 0}) \ =  & \  \ \frac{ p^{2g-1}}{2^{2g-1}} \cdot  \sum_{\theta =1}^{p-1}\frac{1}{\mr{sin}^{2g-2}(\frac{\pi \cdot  \theta}{p})}  \\
 \Big(  = & \ \ \frac{(-1)^{g-1}\cdot p^{2g-1}}{2} \cdot  \sum_{\zeta^p =1, \zeta \neq 1}\frac{\zeta^{g-1}}{(\zeta-1)^{2g-2}} 
 \ \Big). \end{split} \]
 \et
\begin{proof}
Suppose that $X$ is an ordinary (cf. the discussion preceding Lemma 4.4)  proper smooth curve over $k$ classified by a $k$-rational point of $\Mgp$ which lies in the complement of the image of $\Dp_{g, \mbF_p} \setminus {^\circledcirc \Dp_{g, \mbF_p}}$  via the natural projection $\Dp_{g, \mbF_p} \migi \Mgp$  (cf. Theorem 3.3; the discussion preceding Theorem 3.3).
Then  it follows from Theorem 3.3, Proposition 4.3, and Lemma 4.4 that  $\mcQ_{F_*(\mcL^\vee)/X^{(1)}/k}^{2, 0}$ is finite and \'{e}tale  over $k$.

Next, we determine the value of $\mr{deg}_k(\mcQ_{F_*(\mcL^\vee)/X^{(1)}/k}^{2, 0})$.
Denote by $W$ the ring of Witt vectors with coefficients in $k$ and
$K$ the fraction field of $W$.
Since 
$\mr{dim}(X^{(1)}) =1$, which implies that $H^2(X_F, \Omega^\vee_{X^{(1)}})=0$,
 it follows from  well-known generalities concerning  deformation theory  that 
  $X^{(1)}$ may be lifted to a smooth proper curve $X^{(1)}_W$ over $W$ of genus $g$.
In a similar vein, the fact that $H^2(X^{(1)}, \mcE nd_{\mcO_{X^{(1)}}} (F_*(\mcL^\vee)))=0$ implies that  $F_*(\mcL^\vee)$ may be lifted to a vector bundle $\mcE$ on $X^{(1)}_W$.

Now let  $\eta$ be a $k$-rational point of  $\mcQ_{F_*(\mcL^\vee)/X^{(1)}/k}^{2, 0}$ classifying
an injective morphism $i : \mcF \migi F_*(\mcL^\vee)$.
The tangent space to $\mcQ_{F_*(\mcL^\vee)/X^{(1)}/k}^{2, 0}$ at $\eta$ may be naturally identified with
the $k$-vector space $\mr{Hom}_{\mcO_{X^{(1)}}}(\mcF, F_*(\mcL^\vee)/i(\mcF))$,
and the obstruction to lifting $\eta$ to any first order thickening of $\mr{Spec}(k)$ is given by an element of $\mr{Ext}^1_{\mcO_{X^{(1)}}}(\mcF, F_*(\mcL^\vee)/i(\mcF))$.
On the other hand, since, as was observed above, $\mcQ_{F_*(\mcL^\vee)/X^{(1)}/k}^{2, 0}$ is \'{e}tale over $\mr{Spec}(k)$,
it holds that $\mr{Hom}_{\mcO_{X^{(1)}}}(\mcF, F_*(\mcL^\vee)/i(\mcF))=0$, and hence  $\mr{Ext}^1_{\mcO_{X^{(1)}}}(\mcF, F_*(\mcL^\vee)/i(\mcF))=0$ by Lemma 5.3  below.
This implies that
 $\eta$ may be  lifted to a $W$-rational point of $\mcQ_{\mcE/X^{(1)}_W/W}^{2, 0}$, 
 and hence that  $\mcQ_{\mcE/X^{(1)}_W/W}^{2, 0}$ is finite and \'etale over $W$ by Lemma 5.3 and the vanishing of $\mr{Hom}_{\mcO_{X^{(1)}}}(\mcF, F_*(\mcL^\vee)/i(\mcF))$.  Now it follows from  a routine argument 
 that  $K$ may be supposed to be a subfield of $\mbC$.
Denote by $X^{(1)}_{\mbC}$  the base-change of $X^{(1)}_{W}$ via the morphism $\mr{Spec}(\mbC) \migi  \mr{Spec}(W)$ induced by the composite embedding $W \migiincl K \migiincl \mbC$, and 
$\mcE_\mbC$ the pull-back of $\mcE$ via the natural morphism $X^{(1)}_{\mbC} \migi X^{(1)}_{W}$. 
Thus, we obtain equalities
\[ \mr{deg}_k(\mcQ^{2, 0}_{F_*(\mcL^\vee)/X_k/k}) = \mr{deg}_W(\mcQ^{2, 0}_{\mcE/X^{(1)}_W/W})=  \mr{deg}_\mbC(\mcQ^{2, 0}_{\mcE_\mbC/X^{(1)}_\mbC/C}). \]
To prove the required formula, we calculate the degree  $\mr{deg}_\mbC(\mcQ^{2, 0}_{\mcE_\mbC/X^{(1)}_\mbC/\mbC})$ by applying Theorem 5.1.

By  ~\cite{SUN}, \S\,2, p.\,431, Theorem 2.2, $F_*(\mcL^\vee)$ is stable.
Since the degree of $\mcE_\mbC$ coincides with the degree of   $F_*(\mcL^\vee)$, 
$\mcE_\mbC$ is a  vector bundle of degree $\mr{deg}(\mcE_\mbC ) = (p-2)(g-1)$ (cf. the proof of Lemma 5.3). 
On the other hand, 
one verifies easily from the definition of stability and the properness of Quot schemes (cf.  ~\cite{FGA}, \S\,5.5, p.\,127, Theorem 5.14) that $\mcE_\mbC$ is a stable vector bundle.
Next, let us observe that $\mcQ^{2, 0}_{\mcE_\mbC/X^{(1)}_\mbC/C}$ is zero-dimensional (cf. the discussion above), which, by the discussion preceding Theorem 5.1, implies that
$s_2(\mcE_\mbC) = 2(p-2)(g-1)$.  Thus, by choosing the deformation $\mcE$ of $F_*(\mcL^\vee)$ appropriately, 
we may assume, without loss of generality, that $\mcE_\mbC$ is sufficiently general in
$^\text{s}\mcN_{X^{(1)}_\mbC}^{p, (p-2)(g-1)}$
 that Theorem 5.1 holds.  
Now we compute (cf. the discussion preceding Theorem 5.1): 
 \[ \begin{split} e_{\mr{max}}(\mcE_\mbC, 2) = & \ \frac{1}{p} \cdot (\mr{deg}_\mbC(\mcE_\mbC ) \cdot 2 - s_2(\mcE_\mbC))  \\
  = & \ \frac{1}{p} \cdot ((p-2)(g-1) \cdot 2 - 2 \cdot (p-2)(g-1)) \\
 = &  \ 0. \end{split} \]
If, moreover, we write $(a,b)$ for the unique pair of integers  such that $\mr{deg}_\mbC(\mcE_\mbC ) = p \cdot a -b$ with $0 \leq b < p$,
 then it follows from the hypothesis $p>2(g-1)$ that $a = g-1$ and $b= 2 (g-1)$.
Thus,
by applying Theorem 5.1 in the case where the data 
\[ ``(C, \mcV, n, d, r, a, b, e_{\mr{max}}(\mcV, r))"\]
 is taken to be
  \[(X^{(1)}_{\mbC}, \mcE_\mbC, p, (g-1)(p-2), 2, g-1, 2(g-1), 0),\]
   we obtain that
\[ \begin{split}  \mr{deg}_\mbC(\mcQ^{2, 0}_{\mcE_\mbC/X^{(1)}_\mbC/\mbC})  
 &= \frac{(-1)^{(2-1)(2(g-1)2-(g-1)2^2)/p}p^{2(g-1)}}{2!} \cdot \sum_{\rho_1, \rho_2}\frac{(\prod_{i=1}^2 \rho_i)^{2(g-1)-g+1}}{\prod_{i \neq j} (\rho_i - \rho_j)^{g-1}} \\
 \end{split} \] 
 \[ \begin{split} 
&= \frac{(-1)^{g-1}\cdot p^{2g-1}}{2} \cdot  \sum_{\zeta^p=1, \zeta \neq 1}\frac{\zeta^{g-1}}{(\zeta-1)^{2g-2}}  \\
&= \frac{ p^{2g-1}}{2^{g}} \cdot  \sum_{\zeta^p=1, \zeta \neq 1}\frac{1}{(1-\frac{\zeta + \zeta^{-1 }}{2})^{g-1}}  \\
&= \frac{ p^{2g-1}}{2^{2g-1}} \cdot  \sum_{\theta =1}^{p-1}\frac{1}{\mr{sin}^{2g-2}(\frac{\pi \cdot \theta}{p})}.\end{split}  \]
This completes the proof of the required equality.
\end{proof}
\vspace{3mm}
The following lemma was used in the proof of Theorem 5.2.

\vspace{3mm}
\ble \leavevmode\\
 \ \ \  Let $k$ be a field of characteristic $p$, $(X/k, \mbL :=(\mcL, \eta_\mcL))$ a spin curve,  and
  $i:\mcF \migi F_*(\mcL^\vee)$ an injective morphism classified by  a $k$-rational point of $\mcQ_{F_*(\mcL^\vee)/X^{(1)}/k}^{2, 0}$. Write $\mcG := F_*(\mcL^\vee) /i(\mcF)$.
 Then 
 $\mcG$ is a vector bundle on $X^{(1)}$,  and
  it holds that
 \[ \mr{dim}_k(\mr{Hom}_{\mcO_{X^{(1)}}}(\mcF, \mcG) )= \mr{dim}_k(\mr{Ext}^1_{\mcO_{X^{(1)}}}(\mcF, \mcG)). \]
 \ele
\begin{proof}
First, we verify that $\mcG$ is  a vector bundle.
Since $F : X \migi  X^{(1)}$ is faithfully flat, it suffices to verify that the pull-back $F^*\mcG$ is a vector bundle on $X$.
Recall (cf. Lemma 4.2 (i)) that 
the composite $F^*\mcF \migi \mcA_\mcL(=F^*F_*(\mcL^\vee))  \migi \mcA_\mcL/\mcA_\mcL^2$ of the pull-back of $i$ with the natural surjection $\mcA_\mcL \migi \mcA_\mcL/\mcA_\mcL^2$ is an isomorphism.
One verifies easily that this implies that 
the natural composite $\mcA_\mcL^2 \migi \mcA_\mcL \migi F^*\mcG$ is an isomorphism, and hence that $F^* \mcG$ is a vector bundle, as desired.

Next we consider the asserted equality.
Since the morphism $F:X \migi X^{(1)}$ is finite, it follows from well-known generalities concerning cohomology that we have an equality of Euler characteristics $\chi(F_*(\mcL^\vee)) = \chi (\mcL^\vee)$.
Thus, it follows from 
the Riemann-Roch theorem that 
\[ \begin{split} \mr{deg}(F_*(\mcL^\vee)) = & \ \chi (F_*(\mcL^\vee)) - \mr{rk}(F_*(\mcL^\vee))(1-g) \\
= & \ \chi (\mcL^\vee) - p(1-g) \\
= & \ (p-2)(g-1), \end{split} \]
and, since $ \mr{rk} (\mcH om_{\mcO_{X^{(1)}}}(\mcF, \mcG)) = 2 (p-2)$, that 
\[  \begin{split} \mr{deg}( \mcH om_{\mcO_{X^{(1)}}}(\mcF, \mcG)) = & \ 2 \cdot \mr{deg}(\mcG) - (p-2) \cdot \mr{deg}(\mcF)  \\
  = & \ 2 \cdot \mr{deg}(F_*(\mcL^\vee)) - 0 \\
  = &  \ 2 (p-2) (g-1). \end{split} \]
Finally,
by applying the Riemann-Roch theorem again, we obtain equalities
\[ \begin{split} & \mr{dim}_k(\mr{Hom}_{\mcO_{X^{(1)}}}(\mcF, \mcG)) - \mr{dim}_k(\mr{Ext}_{\mcO_{X^{(1)}}}^1(\mcF, \mcG)) \\
 = &  \ \mr{deg}(\mcH om_{\mcO_{X^{(1)}}} (\mcF, \mcG)) + \mr{rk} (\mcH om_{\mcO_{X^{(1)}}} (\mcF, \mcG)) (1-g) \\
 =  & \ 2(p-2)(g-1) + 2(p-2)  (1-g) \\
 =  & \ 0. \end{split} \]
\end{proof}

 Thus, we  conclude the main result of the present paper.

\vspace{3mm}
\bco \leavevmode\\
 \ \ \   Suppose that $p>2(g-1)$.
  Then the degree $\mr{deg}_{\Mgp}(\Dp_{g, \mbF_p})$
 of $\Dp_{g, \mbF_p}$ over $\Mgp$
 is given by the following formula:
 \[ \begin{split} \mr{deg}_{\Mgp}(\Dp_{g, \mbF_p})
 = \  & \  \frac{p^{g-1}}{2^{2g-1}} \cdot  \sum_{\theta =1}^{p-1}\frac{1}{\mr{sin}^{2g-2}(\frac{\pi \cdot \theta}{p})}  \\
 \Big(  \ =  \ & \  \frac{(-1)^{g-1}\cdot p^{g-1}}{2} \cdot  \sum_{\zeta^p =1, \zeta \neq 1}\frac{\zeta^{g-1}}{(\zeta-1)^{2g-2}}
 \ \Big). \end{split} \]

 \eco
\begin{proof}
Let us fix  a spin  curve $(X/k, \mbL)$ for which Theorem 5.2 holds.
Then it follows from  Theorem 5.2 and the discussion at the beginning of \S\,5 that
 \[  \begin{split} \mr{deg}_{\Mgp}(\Dp_{g, \mbF_p}) \ = & \ \ \frac{1}{p^g} \cdot   \mr{deg}_\mbC(\mcQ_{F_*(\mcL^\vee)/X^{(1)}/k}^{2, 0}) \\
 \ = & \ \ \frac{ p^{g-1}}{2^{2g-1}} \cdot  \sum_{\theta =1}^{p-1}\frac{1}{\mr{sin}^{2g-2}(\frac{\pi \cdot \theta}{p})}  \\
 \Big(  \ = & \ \ \frac{(-1)^{g-1}\cdot p^{g-1}}{2} \cdot  \sum_{\zeta^p =1, \zeta \neq 1}\frac{\zeta^{g-1}}{(\zeta-1)^{2g-2}} 
\  \Big). \end{split}\]
\end{proof}

\vspace{6mm}
\section{Relation with other results} \vspace{3mm}

Finally, we discuss some topics related to the main result of the present paper.

\subsection{}

Let $k$ be an algebraically closed field of characteristic $p$ and
$X$  a proper smooth curve over $k$ of genus $g$ with $p>2(g-1)$.
Denote by $F:X\migi X^{(1)}$ the relative Frobenius morphism.
Let $\mcE$ be 
an indecomposable  vector bundle on $X$ of rank $2$ and degree $0$.
If $\mcE$ admits a rank one subbundle of positive degree, then it follows from the definition of semistability that $\mcE$ is not semistable.
On the other hand,  since $\mcE$ is indecomposable, a computation of suitable $\mr{Ext}^1$ groups via Serre duality shows  that the degree of any rank one subbundle of $\mcE$ is at most $g-1$.
We shall say that $\mcE$ is {\it maximally unstable} if
 $\mcE$ admits a rank one subbundle of degree $g-1 (>0)$.
Let us denote by $B$ the set of isomorphism classes of rank $2$ semistable  bundles $\mcV$ on $X^{(1)}$  such that $\mr{det}(\mcV) \cong \mcO_X$, and  $F^*\mcV$ is indecomposable and maximally unstable. 
Then it is well-known (cf., e.g., ~\cite{O1}, \S\,4, p.\,110, Proposition 4.2) that
there is  a natural $2^{2g}$-to-$1$ correspondence between $B$ and the set of isomorphism classes of dormant indigenous bundles on $X/k$.
Thus, Corollary 5.4 of the present paper enables us to calculate the cardinality of $B$, i.e., to conclude that  \[\sharp  B = 2 \cdot p^{g-1}\cdot  \sum_{\theta =1}^{p-1}\frac{1}{\mr{sin}^{2g-2}(\frac{\pi \cdot \theta}{p})}.\]
In the case where $g=2$, this result is consistent with the result obtained in ~\cite{LP}, Introduction, p.\,180, Theorem 2. 

\subsection{}

F. Liu and B. Osserman have shown (cf. ~\cite{LO}, \S\,2, p.\,127, Theorem 2.1) that 
the value $\mr{deg}_{\Mgp}(\Dp_{g, \mbF_p})$ may be expressed as a polynomial with respect to the characteristic $p$ of degree $3g-3$ (e.g., $\mr{deg}_{\mcM_{2, \mbF_p}}(\Dp_{2, \mbF_p})= \frac{1}{24} \cdot (p^3 -p)$, as referred to in Introduction).
In fact, this result may also be obtained as a consequence of Corollary 5.4.
This may not be apparent at first glance, but nevertheless may be verified by applying
either of the following two different (but, closely related) arguments.
\begin{itemize}
\item[(1)]
Let $C$ be a connected compact Riemann surface of genus $g>1$.
Then it is known that the moduli space of S-equivalence classes  (cf. ~\cite{HL}, \S\,1.5, p.\,24, Definition 1.5.3) of rank $2$  semistable bundles on $C$ with trivial determinant
\[{^{\text{ss}}\N{C}{2, \mcO}}\] 
may be represented by a projective algebraic variety  of dimension $3g-3$ (cf.  ~\cite{S}, \S\,8, p.\,333, Theorem 8.1; ~\cite{BEA}, \S\,1, p.\,18; ~\cite{NR}, Introduction),
and that
$\mr{Pic}({^{\text{ss}}\N{C}{2, \mcO}}) $ $\cong \mbZ \cdot [\mcL]$ for a certain ample line bundle $\mcL$ (cf. ~\cite{DN}, \S\,0, p.\,55, Theorem B; ~\cite{BEA}, \S\,2, p.\,19, Theorem 1; ~\cite{BEA}, p.\,21, the discussion at the beginning of \S\,4).
The Verlinde formula,
introduced in ~\cite{V} and proved, e.g., in ~\cite{Fa}, \S\,4, p.\,367, Theorem 4.2,
implies that, for $k = 0, 1, \cdots$, we have an equality
\[  \mr{dim}_\mbC(H^0({^{\text{ss}}\N{C}{2, \mcO}}, \mcL^{\otimes k})) =  \frac{(k+2)^{g-1}}{2^{g-1}} \cdot  \sum_{\theta =1}^{k+1}\frac{1}{\mr{sin}^{2g-2}(\frac{\pi \cdot \theta}{k+2})}  \] 
(cf. ~\cite{BEA}, \S\,5, p.\,24, Corollary).
Thus, for sufficiently large $k$, the value at $k$ of the Hilbert polynomial $\mr{Hilb}_{\mcL}(t) \in \mbQ[t]$ of $\mcL$
coincides with the RHS of the above equality.
On the other hand, it follows from Corollary 5.4 that for an odd  prime $p$, the value at $k=p-2$ of this RHS divided by $2^g$ coincides with the value $\mr{deg}_{\Mgp}(\Dp_{g, \mbF_p})$.
Thus,  the value $\mr{deg}_{\Mgp}(\Dp_{g, \mbF_p})$ (for sufficiently large $p$) may be expressed as $\mr{Hilb}_{\mcL}(p-2)$ for a suitable polynomial
$\mr{Hilb}_{\mcL}(t) \in \mbQ[t]$
of degree $3g-3\ \ \  (=\mr{dim}({^{\text{ss}}\N{C}{2, \mcO}}))$ 
.

\item[(2)]
By comparison to the discussion of (1),
the approach of the following discussion yields a more concrete expression for $\mr{deg}_{\Mgp}(\Dp_{g, \mbF_p})$.
For a pair of positive integers $(n,k)$, we set
\[  V(n,k) := \sum_{\theta =1}^{k-1}\frac{1}{\mr{sin}^{2n}(\frac{\pi  \cdot \theta}{k})}.  \]
Then it follows from ~\cite{Z}, p.\,449, Theorem 1 (i), (ii);  ~\cite{Z}, p.\,449, the proof of Theorem 1 (iii), that
\[     V(n,k) = - \mr{Res}_{x=0} \Big[ \frac{k \cdot \mr{cot}(kx)}{\mr{sin}^{2n}(x)}dx\Big] , \]
where 
$\mr{Res}_{x=0}(f)$ denotes the residue of $f$ at $x=0$.
Thus, $V(n,k)$
 may be computed by considering 
 the relation
 $\frac{1}{\mr{sin}^2(x)} = 1 + \mr{cot}^2(x)$
 and  the coefficient of the Laurent expansion (cf. ~\cite{Z}, p.\,449, the proof of Theorem 1 (iii))
\[ \mr{cot}(x) = \frac{1}{x}+\sum_{j=1}^{\infty}\frac{(-1)^{j}2^{2j}B_{2j}}{(2j)!}x^{2j-1}\]
where $B_{2j}$ denotes the $(2j)$-th Bernoulli number, i.e.,
\[  \frac{w}{e^w-1} = 1 - \frac{w}{2} + \sum_{j=1}^\infty \frac{B_{2j}}{(2j)!}w^{2j}. \] 
In particular, it follows from  an explicit computation that $V(n,k)$ may be expressed as a polynomial of degree $2n$ with respect to $k$.
Thus, the value $\mr{deg}_{\Mgp}(\Dp_{g, \mbF_p})$ ($=\frac{p^{g-1}}{2^{2g-1}} \cdot V(g-1, p)$ by Corollary 5.4) may be expressed as a  polynomial with respect to $p$ of degree $2(g-1)+(g-1)=3g-3$.
Moreover, by applying the above discussion to our calculations,  we obtain the following explicit expressions for the polynomials under consideration:

\[ \begin{split} 
\mr{deg}_{\mcM_{2, \mbF_p}}(\Dp_{2, \mbF_p}) = & \frac{1}{24} \cdot (p^3-p),\\ 
\mr{deg}_{\mcM_{3, \mbF_p}}(\Dp_{3, \mbF_p}) = & \frac{1}{1440} \cdot (p^6+10p^4-11p^2), \\
\mr{deg}_{\mcM_{4, \mbF_p}}(\Dp_{4, \mbF_p})  = & \frac{1}{120960}\cdot (2p^9+21p^7+168p^5-191p^3),\\
 \mr{deg}_{\mcM_{5, \mbF_p}}(\Dp_{5, \mbF_p}) = & \frac{1}{7257600} \cdot (3p^{12}+40p^{10}+294p^8+2160p^6-2497p^4),\\
\mr{deg}_{\mcM_{6, \mbF_p}}(\Dp_{6, \mbF_p}) = & \frac{1}{2048} \cdot \Big(\frac{2}{93555} p^{15}+ \frac{1}{2835}p^{13} + \frac{26}{8505} p^{11} + \frac{164}{8505} p^9 \\
& + \frac{128}{945} p^7 - \frac{14797}{93555} p^5\Big), \\ 
\mr{deg}_{\mcM_{7, \mbF_p}}(\Dp_{7, \mbF_p}) = & \frac{1}{8192} \cdot \Big(\frac{1382}{638512875}p^{18}+ \frac{4}{93555} p^{16} +\frac{31}{70875} p^{14} \\
& + \frac{556}{178605}p^{12} + \frac{3832}{212625}p^{10} +\frac{256}{2079}p^8 -\frac{92427157}{638512875}p^6\Big) \\  
\mr{deg}_{\mcM_{8, \mbF_p}}(\Dp_{8, \mbF_p}) = & \frac{1}{32768} p^7\cdot  \Big( \frac{4}{18243225}p^{14} + \frac{1382}{273648375}p^{12} 
+ \frac{4}{66825}p^{10} \\
& +  \frac{311}{637875}p^8 + \frac{1184}{382725}p^6 +\frac{1888}{111375}p^4 + 
\frac{1024}{9009}p^2 \\
&- \frac{36740617}{273648375} \Big) \\ %
\end{split}\] 
\[ \begin{split}%
 \mr{deg}_{\mcM_{9, \mbF_p}}(\Dp_{9, \mbF_p}) = 
 & 
 \frac{1}{131072}p^8 \cdot \Big(\frac{3617}{162820783125}p^{16} + \frac{32}{54729675}p^{14} 
 \\
 & + \frac{226648}{28733079375}p^{12}  + \frac{2144}{29469825}p^{10} + 
 \frac{4946}{9568125}p^8
 \\
 &+\frac{268864}{88409475}p^6 + \frac{17067584}{1064188125}p^4  + \frac{2048}{19305}p^2 
 \\
 &-\frac{61430943169}{488462349375}
 \Big)
 \\ %
\mr{deg}_{\mcM_{10, \mbF_p}}(\Dp_{10, \mbF_p}) = & \frac{1}{524288}p^{9} \cdot \Big(\frac{87734}{38979295480125}p^{18} + \frac{3617}{54273594375}p^{16}\\
& + \frac{92}{91216125}p^{14} + \frac{2092348}{201131555625}p^{12} + \frac{4042}{49116375}p^{10} \\
& + \frac{18716}{35083125}p^8 + \frac{119654944}{40226311125}p^6 + \frac{16229632}{1064188125}p^4\\
& + \frac{32768}{328185}p^2 - \frac{23133945892303}{194896477400625}\Big).
\end{split}\] 
\end{itemize}

\end{document}